\documentclass[12pt]{article}
\usepackage{amsmath} \usepackage{amsfonts} 

\usepackage{tikz}

\begin{document}

\title{The square of a planar
cubic graph is $7$-colorable}

\author{Carsten Thomassen \thanks{Research partly supported by ERC Advanced Grant GRACOL}\\Department of Applied Mathematics and Computer Science,
\\Technical University of Denmark, DK-2800 Lyngby, Denmark}
\maketitle

\begin{abstract}

We prove the conjecture made by G.Wegner in 1977 that
the square of every planar,
cubic graph is $7$-colorable. Here, $7$ cannot be replaced by $6$.

\vspace{7 mm}


\end{abstract}
\newtheorem{theorem}{Theorem}
\newtheorem{proposition}{Proposition }
\newtheorem{lemma}{Lemma}
\newtheorem{corollary}{Corollary}
\newtheorem{claim}{Claim}
\newtheorem{conjecture}{Conjecture}
\newtheorem{problem}{Problem}
\newtheorem{case}{Case}

\section{Introduction}

We prove the conjecture made by G.Wegner~\cite{w} in 1977,
mentioned by Gionfriddo~\cite{g}
and listed in the monograph by Jensen and Toft~\cite{jt},
that the square of every
planar, cubic graph is $7$-colorable. To see that this bound is best possible, consider first the cubic prism graph with six vertices. Then subdivide an edge which is not contained in a triangle. The square of this graph is a complete graph with seven vertices. Now we take two copies of this graph and add an edge between them so that we obtain a cubic graph.
This cubic graph is planar and its square has chromatic number $7$.

The proof is based on a decomposition method: We color the vertices of the planar, cubic graph by two colors, red and blue, such that the blue square-graph is $3$-colorable, and the red square-graph is planar and hence $4$-colorable, by the $4$-color theorem.

Wegner's $7$-color conjecture proved in the present paper is part of a more general problem on the chromatic number of squares of planar graphs. After submission of the present paper
a number of papers have been written on this subject, see e.g. \cite{bz}, \cite{cj}, \cite{dkns}, \cite{kp}, \cite{ms} and the references in these papers. A computer aided proof of the $7$-color theorem has recently been obtained in \cite{hjt}.

\section{Terminology and notation}
The terminology is the same as in \cite{jt} and \cite{mt}.

A \emph{$k$-path} is a path with $k$ vertices. A
\emph{$k$-cycle} is defined analogously.

In a plane embedding of a connected graph every face boundary is a walk called a \emph{facial walk}.
A \emph{facial path} is a path which is a subgraph of a facial
walk. If $C$ is a cycle in a plane graph, then \emph{the interior of $C$}, denoted $int(C)$, consists of the edges and vertices inside $C$. Thus, an edge joining a vertex in $C$ with a vertex inside $C$ is also in $int(C)$. Sometimes $int(C)$ also refers to a graph, namely the subgraph of $G$ induced by the vertices inside $C$. The precise meaning will always be clear from the context.

If $G$ is a graph, then the
\emph {square $G^2$} of $G$ is obtained from $G$ by adding all edges joining
vertices of distance $2$ in $G$. If we color the vertices of $G$ red
or blue, then the  \emph{red subgraph} (or just the \emph{red graph})   is the subgraph of $G$
induced by the red vertices. The \emph{red square-subgraph} (or just the \emph{red square-graph})  is the
subgraph of $G^2$ induced by the red vertices. Similar notation is used for the blue vertices.

If some vertices of $G$ are colored $1,2,3$ such that the coloring is proper in $G^2$, then we say that \emph{vertex $v$ can see color $i$} if there is a vertex $u$ of color $i$ such that $u$ is a neighbor of $v$ in $G^2$. A \emph{Kempe chain with colors $i,j$} is a connected component in the subgraph of $G^2$ induced by the vertices of colors $i,j$.

We shall also use the following notation: If we have already named a
sequence $v_1,v_2, \ldots$ of vertices in the cubic graph and that
sequence includes say two neighbors of $v_1$, then the neighbor of
$v_1$ which is not in the list is called  \emph{the third neighbor of
$v_1$}. If precisely one neighbor of $v_2$ is in the list, then the
two neighbors of $v_2$ not in the list are called \emph{the two other neighbors of
$v_2$}.

\section{Decomposing the vertex set of a cubic graph}

We shall now indicate the idea in the proof of Wegner's conjecture. We
begin with a conjecture.

\begin{conjecture}
\label{c1}

If $G$ is a $3$-connected, cubic graph, then the vertices of $G$ can be colored blue
and red such that the blue subgraph has maximum degree $1$ (that is,
it consists of a matching and some isolated vertices) and the red
subgraph has minimum degree at least $1$ and contains no $4$-path.

\end{conjecture}

Suppose Conjecture \ref{c1} is true. Assume now that $G$ is planar cubic and
$3$-connected, and consider the blue vertices. Using the fact that the red subgraph has
no isolated vertex, it follows easily that the blue square-graph has
maximum degree at most $3$. As $G$ is $3$-connected it is easy to see
that the blue square-graph contains no complete graph with $4$
vertices. Hence it is $3$-colorable, by Brooks' theorem.
Consider next the red subgraph. Add all edges in the red square-subgraph. Every new edge can be added such that it does not cross any
edge of $G$. Two new edges may cross, though. But, as there is no red facial $4$-path, no two additional edges
cross. So the red square-subgraph is planar and hence $4$-colorable, by
the $4$-Color Theorem. This implies that $G^2$ is $7$-colorable.

The method of this paper is to prove a technical and less elegant version of Conjecture \ref{c1} which is
strong enough, though, to prove Wegner's conjecture.

\section{Decomposing the vertex set of a planar cubic graph}

We shall consider a graph where some vertices are colored red or blue, and some vertices are uncolored. A \emph{forbidden cycle} is a cycle $C'$
such that either the length of $C'$ is not a multiple of $3$ and all the vertices of $C'$ are all blue, or the length of $C'$ is congruent to $2$ modulo $3$ and all but precisely one of the vertices are blue.
Note that the blue vertices cannot be properly colored in three colors in the square of a forbidden cycle.  We say that a cycle is a \emph{dangerous cycle} if and only if it has at least one non-blue vertex, and if we change the color of any vertex from non-blue to blue, then we obtain a forbidden cycle.

It is easy to see that a cycle $C''$ is dangerous if and only if $C''$ has length congruent to $1$ modulo $3$ such that
each vertex of $C''$, except precisely one, is blue, or $C''$ has length congruent to $2$ modulo $3$ such that
each vertex of $C''$, except precisely two, is blue.

\begin{theorem}
\label{t1}

Let $G$ be a connected planar graph with chordless outer cycle $C$. Assume that the following hold:

$(c_1)$: If $v$ is a vertex in $int(C)$ and $E$ is a set of at most two edges in $int(C)$, then $G-E$ has a path from $v$ to $C$.

$(c_2)$: All vertices of $G$ have degree at most $3$, and all vertices in $int(C)$ have degree precisely $3$.

$(c_3)$: Each vertex of $C$ is precolored by red or blue such that at most one vertex $b_0$ of $C$ is blue.

$(c_4)$: If $b_0$ exists, then at least one of its neighbors on $C$ has degree $2$ in $G$.

$(c_5)$: If $b_0$ does not exist, that is, all vertices of $C$ are red, then some vertex $r_0$ on $C$ is called either left-forbidden or right-forbidden or $4$-forbidden.

$(c_6)$: $G-V(C)$ is connected and contains a vertex joined to $b_0$ (if $b_0$ exists) or $r_0$ (otherwise).

$(c_7)$: Every vertex $v$ in $int(C)$ which has a red neighbor on $C$ distinct from $r_0$ is colored blue. If $r_0$ exists and is either right-forbidden or left-forbidden and its neighbor $r_0'$ in $int(C)$ has two neighbors in $int(C)$, then one of these neighbors is precolored blue as follows: If $r_0$ is right-forbidden (respectively left-forbidden) and the path $r_0r_0'a$ turns sharp right (respectively left) at $r_0'$, then $a$ is precolored blue. No other vertex in $int(C)$ is precolored.
(To clarify, we emphasize some special consequences of this:
If the neighbor $r_0'$ of $r_0$ in $int(C)$ has two neighbors on $C$, then $r_0'$ is precolored blue, and the third neighbor of $r_0'$
is not precolored unless that neighbor is joined to $C-r_0$. If $r_0$ is $4$-forbidden and $r_0'$ has two neighbors in $int(C)$, then $r_0'$
is not precolored. Also, a neighbor of $r_0'$ in $int(C)$
is not precolored unless that neighbor of $r_0'$ has a neighbor on $C-r_0$.)

$(c_8)$: There is no forbidden cycle, and there is no dangerous cycle in $G$, except possibly a dangerous cycle that contains $b_0$ and precisely one other vertex of $C$.

$(c_9)$: If $r_0$ exists and is right-forbidden or left-forbidden and $r_0'$ is the neighbor of $r_0$ in $int(C)$ and $r_0'$ has no other neighbor in $C$, then $G$ has a facial cycle which contains $r_0'$ and is disjoint from $C$.

\vspace{3mm}

Then the red-blue coloring (of $V(C)$ and the blue vertices inside $C$) can be extended to a red-blue coloring of $V(G)$ such that the following conditions hold:

 \vspace{3mm}

$(i)$: $G$ has no red facial $4$-path whose edges are in $int(C)$. (Note that, if $r_0$ exists and is left-forbidden, then $G$ has no red facial $3$-path starting at $r_0$ and with edges in $int(C)$ and turning sharp left at the neighbor of $r_0$. In other words, we allow a red facial $3$-path starting at $r_0$ and with edges in $int(C)$ provided it turns sharp right at the neighbor of $r_0$. Similarly, if $r_0$ exists and is right-forbidden, then $G$ has no red facial $3$-path starting at $r_0$ and with edges in $int(C)$ and turning sharp right at the neighbor of $r_0$.
 If $r_0$ is neither right-forbidden nor left-forbidden, then the only condition on $r_0$ is that there is no red facial $4$-path starting at $r_0$, and therefore we call it $4$-forbidden.)

$(ii)$: the blue vertices of $G$ can be $3$-colored with colors $1,2,3$ such that this coloring is proper in $G^2$.

$(iii)$: If $r_0$ is right-forbidden or left-forbidden, and the neighbor $r_0'$ of $r_0$ in $int(C)$ has no other neighbor on $C$, then the red-blue coloring can be chosen such that $r_0'$ is red.

\end{theorem}

Note that $(c_1),(c_2),(c_6)$ imply that $G$ is obtained from a planar, cubic, $3$-connected graph by subdividing edges on the outer cycle. The condition $(c_1)$ is included only to reduce the amount of case analysis. It is a triviality to check $(c_1)$ in the induction steps. We may assume that each edge in the outer cycle, except possibly one, is subdivided many times, since subdividing edges of $C$ (except possibly one) does not affect the conditions nor conclusion of Theorem \ref{t1}. This is useful when we wish to add an edge from $int(C)$ to $C$. If $b_0$ exists, then, by $(c_4)$, $b_0$ may be joined to precisely one vertex $d$ on $C$ which has degree $3$ in $G$. We are not allowed to subdivide the edge $b_0d$. The reason for the condition $(c_9)$ is that we keep $int(C)$ connected if we delete $r_0'$. The reason for $(iii)$ is that it is convenient when we color the blue vertices by colors $1,2,3$ in the induction step.

In the induction step we sometimes introduce a new blue vertex. This will never create a forbidden cycle, but it may create a new dangerous cycle which we then have to dispose of. In case we make more than one vertex blue we make sure that the other new blue vertices are not part of a dangerous cycle. This is done by letting the other new blue vertices have two neighbors on the outer cycle.

\vspace{3mm}

\emph{Proof of Theorem \ref{t1}}

The proof is by induction on the number of edges inside $C$. Suppose (reductio ad absurdum) that Theorem \ref{t1} is false. Select a counterexample such that the number of edges in $int(C)$ is minimum.

%The idea in the proof is to first extend $C$ to a subgraph $H$, then color $H$, and finally apply the induction hypothesis to each component of $G-V(H)$. Those vertices in $H$ which are blue because of $(c_7)$ are denoted $b_j$. The other vertices in $H$ are denoted $r_j$. Almost all of the vertices $r_j$ are red when we apply induction. After the induction a few of them may change color from red to blue. To avoid confusion, we call the final red-blue coloring of $V(H)$ the \emph{post-coloring} whereas we denote the red-blue coloring of $V(H)$ during the induction step as the \emph{pre-coloring}.

\vspace{3mm}

We may assume the following:

\vspace{3mm}

 { \bf Claim 1}: $int(C)$ does not contain a path $v_1v_2v_3v_4v_5$ such that each of $v_1,v_2,v_3,v_4,v_5$ is joined to $C-r_0$ (if $r_0$ exists) or $C-b_0$ (if $b_0$ exists).

 \vspace{3mm}

 For, if such a path exists, then we delete the vertices $v_2,v_3,v_4$ add the edge $v_1v_5$ and use induction. If
   $v_1,v_5$
 get blue colors $1,2$, respectively, then $v_1,v_2,v_3,v_4,v_5$ can be colored $1,2,3,1,2$, respectively. This contradiction proves Claim 1.

\vspace{3mm}

{\bf Case $1$: $r_0$ exists and is either right-forbidden or left-forbidden.}

Assume without loss of generality that $r_0$ is right-forbidden. Let $r_0'$ be the neighbor of $r_0$ in $int(C)$.

\vspace{3mm}

{\bf Subcase $1.1$: $r_0'$ has at least two neighbors on $C$.}

\vspace{3mm}

If $r_0'$ has three neighbors on $C$, then
  $(c_6)$ implies that
$int(C)$ consists of $r_0'$ only and there is nothing to prove. So assume that $r_0'$ has a neighbor $a$ in $int(C)$. We define a new graph $G'$ by deleting $r_0'$ and adding the edge $ar_0$. We now apply the induction hypothesis to $G'$ where $r_0$ plays the role of $b_0$. We subdivide some edges of $C$ incident with $r_0$ so that there is no dangerous cycle in $G'$ and such that $(c_4)$ is satisfied. After the induction we transfer the blue color of $r_0$ to $r_0'$, and we give $r_0$ its red color back.

\vspace{3mm}

{\bf Subcase $1.2$: $r_0'$ has only one neighbor $r_0$ on $C$.}

\vspace{3mm}

Let $a,b$ be the neighbors of $r_0'$ in $int(C)$ (and thus distinct from $r_0$) such that $b$ is blue, that is, the path $r_0r_0'b$ turns sharp right at $r_0'$.
 (In other words, the facial path $r_0r_0'b$ is part of a facial cycle traversed clockwise.)
 Let $c,d$ be the neighbors of $r_0$ on $C$ such that the path $cr_0d$ is anticlockwise around $C$.
 Hence the paths $cr_0r_0'a$ and $dr_0r_0'b$ are facial.
 Now we define a new graph $G'$ with a new outer cycle $C'$ as follows: We delete $r_0,r_0'$,
   we add two new vertices $x,y$, we add a path $cxyd$, and we add the edges $xa,yb$.
  Let $u,v$ be the neighbors of $a$ distinct from $x$ such that $cxau$ is a facial path.

\vspace{3mm}

{\bf Subcase $1.2.1$: $u$ is in $C$.}

\vspace{3mm}

Let $G''$ be obtained from $G'$ by deleting $a$ and adding the edge $xv$. We apply induction to $G''$ where $x$ plays the role of $b_0$. (We explain below why induction is possible.)
%Before the induction we need to argue why we do not create a dangerous or forbidden cycle when $x$ is made blue. The reason is that such a cycle would
%contain the edges $vx,xy$. But the corresponding in $G$ has the same length and the same number of blue vertices, a contradiction.
After the induction we transfer the blue color of $x$ (which is distinct from the blue color of $b$) to $a$, and we let $r_0'$ be red. Now we argue why we can apply induction to $G''$. The only problem is that there may be a forbidden or dangerous cycle $C''$ in $G''$. As $G$ has no forbidden or dangerous cycle, $C''$ must contain $x$ or $y$ or both. As we may subdivide the edges $xc,yd$, we may assume that $C''$ contains the path $vxyb$. The corresponding cycle in $G$ contains the path $var_0'b$ which has the same length and same number of blue vertices as $C''$. Hence that cycle in $G$ is forbidden or dangerous, a contradiction to $(c_8)$.

So, we may assume that

\vspace{3mm}

{\bf Subcase $1.2.2$: $u$ is not in $C$.}

\vspace{3mm}

Note that also $v$ is not in $C$ because of $(c_9)$. Possibly, $v,b$ are neighbors. Possibly $v=b$. By $(c_1)$, $u \neq b$, and $u,b$ are nonneighbors.
%because $G$ is a subdivision of a $3$-connected graph.

\vspace{3mm}

We now try to apply induction to $G'$ where $x$ plays the role of $r_0$ and is left-forbidden. If
 $G'$ satisfies $(c_1)-(c_9)$, then we apply induction. In that case $a$ will be red by (iii), and then we have completed Case 1 because we can use the coloring of $G'$
 (which satisfies (i),(ii),(iii)) for $G$.

 It is a simple matter to show that $G'$ satisfies $(c_1)-(c_7)$. We now consider the cases where
 $(c_9)$ or $(c_8)$ fail for $G'$. These are the cases $1.2.2.1,1.2.2.2$, respectively, below.

\vspace{3mm}

{\bf Subcase $1.2.2.1$:
 The facial cycle of $G'$ which contains $a$ but not $x$ intersects the outer cycle $C'$.}

\vspace{3mm}

In this case $a$ is a cutvertex in $int(C')=G'-V(C')$. Let $H_1,H_2$ be the two components of $int(C')-a$ such that $H_1$ contains $u$ and hence $H_2$ contains $b,v$. We now apply induction to $H_1,H_2$ separately. First we draw $H_i$ inside a cycle $C_i$ for $i=1,2$. The vertices in $H_i$ which have a neighbor on $C'$ each has a neighbor on $C_i$. Moreover $u$ is joined to a vertex $a_1$ in $C_1$ and $v$ is joined to a vertex $a_2$ in $C_2$. We may draw $H_2,C_2$ such that the neighbor of $b$ in $C_2$ is a neighbor of $a_2$. Now we use induction such that $a_i$ plays the role of $b_0$ for $i=1,2$. (It is easy to see that induction is possible. The only problem is a possible forbidden or dangerous cycle. There is no forbidden cycle as that would be dangerous in $G$. There may indeed be a dangerous cycle containing a path of length $2$ from $b$ to $a_2$. Such a dangerous cycle is allowed.) We may assume that $a_1,a_2$ have the same blue color, and we give that color to $a$ in $G$. All other colors are transferred to $G$ in the obvious way. Note that
$a,b$ get distinct blue colors. It is possible that $u,v$ have the same blue color. In that case we interchange the two blue colors (distinct from the color of $a_1$) in $H_1$ so that $u,v$ get distinct blue colors.

\vspace{3mm}

  As noted above, $G'$ satisfies$(c_1)-(c_7)$. As have disposed of Subcase 1.2.2.1 we may assume that $G'$ also satisfies $(c_9)$. Now we consider the subcase where $G'$ does not satisfy $(c_8)$.

 \vspace{3mm}

{\bf Subcase $1.2.2.2$: The facial cycle of $G'$ which contains $a$ but not $x$ does not intersect the outer cycle $C'$.}

\vspace{3mm}

As $G'$ does not satisfy $(c_8)$, $G'$ contains a dangerous cycle $C_d$.
  Clearly, $C_d$ contains $u$, since $G$ has no dangerous cycle.
 %By $(c_1)$, $C_d$ does not contain $b$.

 \vspace{3mm}

{\bf Subcase $1.2.2.2.1$: $C_d$ contains $a$.}

\vspace{3mm}

As $C_d$ becomes dangerous when we make $u$ blue
it follows that $u$ is not joined to $C'$ (since otherwise, $C_d$ would be dangerous in $G$.) Also, $v$ is not joined to $C'$ (unless $v=b$) because each of $G$ and $G'$ satisfies $(c_9)$. It follows that all vertices of $C_d$ except $v,u$ are joined to $C'$ (unless $v=b$). We may assume that $u,v$ are not neighbors since otherwise,
 $C_d$ has length $2$ modulo $3$ and hence, by Claim 1,
$C_d$ is of the form $uavbz_1u$ or $uavbz_1z_2z_3z_4u$, where $z_1,z_2,z_3,z_4$ are joined to $C'$. In the former case the edge from $z_1$ to $C$ and the edge $r_0r_0'$ violate $(c_1)$.
Therefore $C_d=uavbz_1z_2z_3z_4u$.
Hence $int(C')$ consists of $C_d$ and the edge $uv$. In this case $r_0'$ is the only vertex in $int(C)$ not in $C_d$. It is easy to complete the proof with $r_0',a,v$ being the only red vertices in $int(C)$. Finally, the neighbors of $u,v$ not in $C_d$ are outside $C_d$ because of $(c_1)$.
\vspace{3mm}

{\bf Subcase $1.2.2.2.1.1$: $v \neq b$.}

\vspace{3mm}

  Because $u,v$ are the only vertices of $C_d$ without neighbors in $C'$, $C_d$ is a facial cycle. Also, $u,v$ are nonconsecutive on $C_d$.
As $v \neq b$, $int(C')-V(C_d)$
has precisely two components $H_1,H_2$ where $H_2$, say, contains $b$ and a neighbor of $v$, and $H_1$ contains a neighbor of $u$. The components $H_1,H_2$ are both outside $C_d$ because $C_d$ is facial. We now try to apply induction to $H_1,H_2$ separately. First we draw $H_i$ inside a cycle $C_i$ for $i=1,2$.
 Each vertex of $H_i$ has as many neighbors on $C_i$ as it has on $C'$.
Moreover the neighbor of $u$ in $H_1$ is joined to a vertex (which we also call $u$) in $C_1$, and the neighbor of $v$ in $H_2$ is joined to a vertex (which we also call $v$) in $C_2$. We let each of $u,v$ play the role of $r_0$, and we call each of $u,v$ $4$-forbidden. After the induction we make $a$ blue and $r_0'$ red.
We then color the blue vertices of $C_d$ by the colors $1,2,3$ starting with $a$ (which gets a color distinct from those of $b$ and the third neighbor of $v$) and then the other neighbor of $v$ finishing with the other neighbor $u'$ of $u$ in $C_d$. When we make this blue coloring we ignore the color of the third neighbor $u''$ of $u$. If $u''$ is blue we permute the blue colors in $H_1$ so that $u''$ gets a blue color distinct from those of $a,u'$.

\vspace{3mm}

  {\bf Subcase $1.2.2.2.1.2$: $v=b$.}
\vspace{3mm}

If $int(C')-V(C_d)$ has precisely two components $H_1,H_2$ where $H_1$ contains a neighbor of $u$, then we repeat the argument in Subcase $1.2.2.2.1.1$ except that now $u$
 plays the role of $b_0$.
The vertex of $C_d$ joined to $H_2$ again plays the role of $r_0$ and is $4$-forbidden. After we have applied the induction hypothesis we color the blue vertices of $C_d$ with colors $1,2,3$ so that there is no blue color conflict with $H_1$. The vertex in $H_2$ joined to $C_d$ may be blue and it may have a color conflict with a blue vertex in $C_d$. In that case we permute the blue colors in $H_2$ so that this conflict disappears.

We now consider the case where $int(C')-V(C_d)$ has precisely one component $H_1$ containing a neighbor of $u$. Consider the case where $C_d: uavw_1w_2u$ where $u$ has a neighbor $u'$ in $H_1$, and $w_2$ has a neighbor $w_2'$ in $H_1$, and $w_1$ is joined to $C$. (It is also possible that $C_d$ has length $4$ in which case $w_2$ does not exist and $u$ is the only vertex of $C_d$ joined to $H_1$ because a dangerous cycle with $4$ vertices has only one vertex that is not precolored blue. That case is similar and easier. It is also possible that $C_d$ has length $7$ in which case $u$ is the only vertex of $C_d$ joined to $H_1$. Again, that case is similar and easier.
%Or $C_d$ may have $7$ vertices with $4$ consecutive vertices joined to the outer cycle.
Note that $C_d$ cannot have length $6$ by the definition of a dangerous cycle,
 and $C_d$ cannot have length at least $8$ because of Claim 1.)
Now we apply induction to $H_1, u$, and an outer cycle $C''$.
 Each vertex of $H_1$ has as many neighbors on $C''$ as it has on $C$.
The vertex $u$ is joined to two vertices of $C''$ which implies that $u$ is precolored blue, and
 $C''$ can be suitably subdivided so that $u$ is
not contained in any dangerous cycle. $w_2'$ is joined to a vertex which is also called $w_2$. This vertex $w_2$ plays the role of $r_0$ and is
 $4$-forbidden.
After the induction $r_0',a$ are made red, and we can give $w_1,v$ blue colors. This coloring satisfies the conclusion of the theorem except that $w_2',u$ may have the same blue color. In that case we first try to interchange colors of $u,a$ so that the color conflict between $w_2',u$ no longer exists. The problem with this color change is that we may create a red facial $4$-path containing $u$. Then we let $u$ keep its blue color, but we recolor $u,w_1,v$ in that order so that we obtain a proper $3$-coloring in the blue square-graph.

%As $C_d$ has length at least $5$ the only problem that may occur with the blue coloring is that we cannot color $a$ because the neighbors of $u,v$ in $H_1,H_2$, respectively, are blue. In that case we make $a$ red instead.

%We only need to argue why we can color $H_1,H_2,u,v$ such that $u,v$ are red. As mentioned above we try to use induction such that each of $u,v$ plays the role of $r_0$. When we use induction, the neighbors of $u,v$ in $H_1,H_2$ become red unless they are joined to the outer cycle. But, we do not need that here. So, we can repeat the argument for $G'$ to each of $H_1,H_2$: Either we complete the proof by induction, or one or both of $H_1,H_2$ contains a cycle with the same property as $C_d$, namely that it becomes dangerous if we make the neighbor of $u$ or $v$ blue.
%By repeating this argument a finite number of times we obtain the desired coloring. This completes Case 1.

\vspace{3mm}

{\bf Subcase $1.2.2.2.2$: $C_d$ does not contain $a$.}

\vspace{3mm}

Let $C_d=uz_1z_2 \ldots z_qu$ traversed anticlockwise. The assumption of Subcase $1.2.2.2$ implies that $z_1$ and its third neighbor $z_1'$ are in $int(C)$. All vertices of $C_d$ except $u,z_1$ and possibly one more, are joined to $C'$ because of the assumption on $C_d$.
\vspace{3mm}

  {\bf Subcase $1.2.2.2.2.1$: $u,v$ are neighbors.}

 \vspace{3mm}

 In this subcase $C_d$ contains the edge $uv$. If $C_d$ does not contain $b$, then,
  because the other neighbor of $v$ in $C_d$ is in the facial cycle through $r_0'$ guaranteed by $(c_9)$,
 $C_d:uvz_2z_3z_4u$ where $z_3,z_4$ are joined to $C$.
 The subgraph of $G$ induced by $C,C_d,a,r_0'$ has precisely one face containing vertices of $G$, by $(c_9)$. We apply induction to that face with $z_2$ being $4$-forbidden. After the induction,
   we change the colors of $a,z_3$ to blue, we color $z_4$ blue as well, and we color $u$ red.
 If the third neighbor of $z_2$ is red, then also $v$ is changed to blue. Otherwise it stays red. Now it is easy to color the blue vertices with the three blue colors.
 If $C_d$ contains $b$, then
   either
 $C_d:uvbz_3z_4z_5z_6u$ where $z_3,z_4,z_5,z_6$ are joined to $C$, and $r_0',a$ are the only vertices in $int(C)$ but not in $C_d$, a configuration which is easy to dispose of, or else $C_d:uvbz_3z_4z_5z_6z_7u$ where four of $z_3,z_4,z_5,z_6,z_7$ are joined to $C$.
   (Again, by $(c_1)$ and Claim 1, there must be either four or five $z_i's$. Also note that $b$ must be the immediate successor of $uv$ on $C_d$ since otherwise we would get a contradiction to $(c_1)$.)
 Then the subgraph of $G$ induced by $C_d,a,r_0'$ has a unique face with vertices of $G$ and the interior of that face has a vertex joined to some $z_i$, $3 \leq i \leq 7$.  We apply induction to that face with $z_i$ playing the role of $r_0$ being $4$-forbidden.
\vspace{3mm}

{\bf Subcase $1.2.2.2.2.2$: $u,v$ are not neighbors, but $C_d$ contains $b$.}

\vspace{3mm}

Then $C_d$ is of the form $uz_1z_2z_3z_4z_5z_6z_7u$ where $z_3=b$, and $z_4,z_5,z_6,z_7$ are joined to $C$.
 (It is not possible that $b$ equals $z_1$ or $z_2$ because of the connectivity condition $(c_1)$. It is not possible that $b$ equals $z_j$ with $j>3$ because $C_d$ is dangerous. Note that $C_d$ cannot have length $4$ or $5$ or $7$ because of $(c_1)$.
Also note that the $5$ consecutive blue vertices on $C_d$ joined to $C'$ does not contradict Claim 1 because only $4$ of them are joined to $C$.) Now we apply induction to the cycle (in $G$) $uar_0'bz_2z_1u$  and its interior, where $a$ plays the role of $b_0$. After the induction we make $b=z_3$ blue, and we color
 $z_3,z_4,z_5,z_6$ (in that order) with colors $1,2,3$.
 So, we may assume that $C_d$ does not contain $b$.

\vspace{3mm}

In order to complete Case 1 only the following subcase remains.

\vspace{3mm}

{\bf Subcase $1.2.2.2.2.3$: $u,v$ are not neighbors, and $C_d$ does not contain $b$.}

\vspace{3mm}

%If $C_d$ contains $v$, then $v=z_1$, {\color{red} because of $(c_1)$}. Then $z_2$ is the only vertex of $C_d$  {\color{red} other than $u,v$} not joined to $C'$. We now consider the subgraph of $G'$ induced by $C',C_d,a$. We apply induction to the face containing $b$, and we let $z_2$ play the role of $b_0$. After the induction we make $z_1,u$ red, and we make $a$ blue. {\color{red} We change the blue color of $z_2$ to red unless this creates a red facial $4$-path.} It is then easy to color the blue vertices by $1,2,3$. So assume that $v$ is not in $C_d$.

As $C_d$ is dangerous when $u$ is made blue, it has at most three vertices not joined to $C'$. One of these is $u$. Another is $z_1$ by the assumption of Subcase $1.2.2.2$. As $u,v$ are not neighbors, $z_1 \neq v$. Now $C_d$ does not contain $v$ because of $(c_1)$.
 Let $z_1'$ be the neighbor of $z_1$ not in $C_d$.

We now form a new graph $G''$ from $G'$ as follows: We first delete $C_d$ and the vertex $a$.
If $int(C')-V(C_d)-a$ is disconnected we focus on the component containing $b$ and ignore
 for
the moment the other component. Then we add an edge from $v$ to $x$, and we add an edge from $z_1'$ to $C'$. If $z_2$ has a neighbor $z_2'$ in the graph under consideration, then we put $z_2$ back and we add two edges from $z_2$ to $C'$. Now we try to use induction with $x$ playing the role of $b_0$. (There may be one more component of $int(C')-V(C_d)-a$ but that component is easy to dispose of as explained below.)

Let us first consider the case where the induction is possible,
that is, $G''$ has no dangerous cycle and hence satisfies $(c_1)-(c_9)$. Then we give $a$ the blue color of $x$ so that $a,b$ have distinct blue colors. We let $u,z_1$ be red. It is possible that $z_2$ has the same blue color as $z_1'$. In that case we try to make $z_2$ red. If that fails it is because we create a red facial $4$-path in which case $z_2$ can see only one blue color except that of $z_1'$. But then $z_2$ can get another blue color so that it has no color conflict with $z_1'$. Then we color the blue vertices of $C_d$ which is possible because we only have to watch the colors of $a$ and $z_2$ or $z_2'$.

If $z_2$ is joined to $C$
(or to a component of $int(C')-V(C_d)-a$ not containing $b$),
and there is a $z_i$ ($i \geq 2$) which is not joined to $C$, then we argue similarly except that we now have to apply induction (with $z_i$ playing the role of $b_0$) to the component of $int(C')-V(C_d)-a$ not containing $b$.

So we may assume that it is not possible to use induction to $G''$. This means that we create a dangerous cycle $C_d'$ in $G''$ when $z_1'$ becomes blue. The cycle $C_d'$ is disjoint from $C_d$. However, $C_d'$ may contain $v$ or $b$ or both.

\vspace{3mm}

{\bf Subcase $1.2.2.2.2.3.1$: The edge $z_1z_1'$ is a bridge in $int(C')-a$.}

\vspace{3mm}

We consider three subcases.

\vspace{3mm}

{\bf Subcase $1.2.2.2.2.3.1.1$:  Each of $z_2,z_q$ is joined to $C$}

\vspace{3mm}

We form $G''$ as in Subcase $1.2.2.2.2.3$ above except that we now delete $a$ and $C_d-z_1$. Again, we add the edge $vx$ and let $x$ play the role of $b_0$. But now we add two edges from $z_1$ to $C$ forcing $z_1$ to be blue. Then we use induction. After the induction we give $a$ the blue color of $x$, and we make $u$ red. There may be a vertex $z_i$
 $(i \geq 3)$
 on $C_d$ such that the edge $z_iz_i'$ from $z_i$ leaving $C_d$ is a bridge in $int(C)$. We then apply induction to the component $Q$ of $int(C)-z_iz_i'$ not containing $a,b$ and with $z_i$ playing the role of $r_0$ and being $4$-forbidden.

 After the induction $z_1,a$ may have the same blue color. We try to make $z_1$ red. If this is possible, it is easy to complete the coloring. So assume that it is not possible, that is, we create a facial $4$-path containing $z_1$. So $z_1'$ and one of
   its
 two other neighbors are red. Hence we can recolor $z_1$ so that it has a blue color that does not conflict with any other blue color. Then we color $z_2, \ldots z_q$ (except $z_i$ which is red) using colors $1,2,3$, ignoring $z_i'$. Then $z_i'$ is the only vertex having color conflicts with other blue vertices. We avoid this by permuting the blue colors in $Q$.

\vspace{3mm}

{\bf Subcase $1.2.2.2.2.3.1.2$: $z_2$ is not joined to $C$}

\vspace{3mm}

Since the vertices $z_3, \ldots ,z_q$ are all joined to $C-r_0$, the path with these vertices has
  at most $4$ vertices, by Claim 1.
 As $C_d$ is dangerous (when $u$ is made blue) it has length $5$, that is, $q=4$. We form $G''$ as in Subcase $1.2.2.2.2.1.1$ above, that is, we delete $a$ and $C_d-z_1$. We add the edge $vx$ and let $x$ play the role of $b_0$. We add two edges from $z_1$ to $C$ forcing $z_1$ to be blue. Then we use induction. After the induction we give $a$ the blue color of $x$, and we make $u$ red. Let $z_2'$ be the neighbor of $z_2$ not in $C_d$. We apply induction to the component $Q$ of $int(C)-z_2z_2'$ not containing $a,b$ and with $z_2$ playing the role of $r_0$ and being $4$-forbidden. After the induction $b,a,z_1,z_3,z_4$ are blue, and $u,r_0',z_2$ are red. The vertices $z_3,z_4$ do not yet have a color $1,2,3$,
   and
  $z_1,a$ may have the same blue color. We try to make $z_1$ red and recolor $z_2,Q$ such that $z_2$ now plays the role of $b_0$. If this is possible, it is easy to complete the coloring. So assume that it is not possible, that is, we create a red facial $4$-path containing $z_1$. So $z_1'$ and one of its two other neighbors are red. Hence we can recolor $z_1$ so that it has a blue color that does not conflict with any other blue color except that of the neighbor $z_2'$ of $z_2$ (which we ignore at the moment). Then we color
  $z_4,z_3$
  using colors $1,2,3$, again ignoring $z_2'$. Then $z_2'$ is the only vertex having color conflicts with others. We avoid this by permuting the blue colors in $Q$.

\vspace{3mm}

{\bf Subcase $1.2.2.2.2.3.1.3$: $z_q$ is not joined to $C$}

\vspace{3mm}

Again, $q=4$. In this case we form a graph $G'''$
  from $G'$
 as follows: We delete $C_d$ and focus on the component of $int(C')-V(C_d)$ containing $a$. We add an additional edge from $a$ to $C$ forcing $a$ to become blue. We also add an edge from $z_1'$ to a vertex of $C$ which we call $z_1$ (with a slight abuse of notation). Then we apply induction to the resulting graph with $z_1$ playing the role of $r_0$ and being $4$-forbidden. We also apply induction to the other component
of $int(C')-V(C_d)$ joined to $z_q=z_4$ such that $z_4$ plays the role of $r_0$ and is $4$-forbidden. Now $a$ is blue and $z_1,z_4$ are red. We also make $u$ blue. Recall that $r_0'$ is red, and $b$ is blue. In order to avoid the possible color conflict between $a,b$ we try to make $a$ red. This is possible unless we create a red facial $4$ path starting at $r_0'$ or $a$. In either case $a$ can see only one blue vertex except $b$ and $u$. So we can give $a$ a blue color so that it has no color conflict (as $u$ does not yet have a color $1,2,3$). We now give the blue vertex $u$ a color $1,2,3$ such that it has no blue color conflict except possibly with the neighbor $z_4'$ of $z_4$ outside $C_d$. This is possible because we only need to watch $a$ or $v$ (but
 not both) and $z_1'$. Then we color the blue vertices $z_2,z_3$ again ignoring the neighbor $z_4'$ of $z_4$ outside $C_d$. That vertex $z_4'$ can be disposed of by permuting the blue colors in the component of $int(C)-z_4z_4'$ not containing $a,b$.

\vspace{3mm}

{\bf Subcase $1.2.2.2.2.3.2$: The edge $z_1z_1'$ is not a bridge in $int(C')-a$.}

\vspace{3mm}

Since $int(C')-a-z_1z_1'$ is connected it has a path from $C_d$ to $C_d'$. Hence $C_d$ has length $2$ modulo $3$,
  and hence length $5$ by Claim 1
 and the path must leave $C_d$ at $z_2$ and enter $C_d'$ at a vertex $w_1$ not joined to $C$.

\vspace{3mm}

{\bf Subcase $1.2.2.2.2.3.2.1$: $b$ is not contained in $C_d'$.}

\vspace{3mm}

Since $int(C)-a$ is connected it has a path from $C_d'$ to $b$. That path must leave $C_d'$ at a vertex not joined to $C$. From these observations it follows that also $C_d'$ has length $2$ modulo $3$. We may assume that $C_d,C_d'$ each has length $5$
  by Claim 1
. (Otherwise we replace an appropriate blue $3$-path by an edge and use induction.) So $C_d=uz_1z_2z_3z_4u$, and $C_d'=z_1'w_1w_2w_3w_4z_1'$
where $w_2,w_3$ are joined to $C$, and the above-mentioned path from $C_d'$ to $b$ leaves $C_d'$ at $w_4$. In particular,
$C_d'$ does not contain $v$.

%We now aim at a coloring where $b,a,z_1$ are blue and $u,z_2$ are red. For this we delete from $G'$ the vertices of $C_d-z_1$. We add an edge from the neighbor $z_2'$ of $z_2$ outside $C_d$ to a vertex on $C$ which we also call $z_2$, and this vertex will play the role of $r_0$ and it will be $4$-forbidden. We add two edges from $z_1$ to $C$ forcing $z_1$ to be blue. We add an additional edge from $a$ to $C$ so that $a$ gets degree $3$. Then we apply induction. If possible we make $a$ red which will make the blue coloring of $z_1,z_3,z_4$ easy because $z_1$ can see at most one blue color on $C_d'$. So assume that $a$ cannot be made red. Then $v$ is red and has at least one red neighbor. Hence, $a$ can see at most one blue neighbor of $v$, and this implies that we can give $a$ another blue color (if necessary) so that $a$ has no color conflict with $b$ or a neighbor of $v$. Now there remains the blue vertices $z_1,z_3,z_4$ which are the only blue vertices that have not yet been colored $1,2,3$. We may assume that the neighbor $z_2'$ of $z_2$ distinct from $z_1,z_3$ is blue and of color $1$, since otherwise, it is easy to color $z_1,z_3,z_4$.

We form a new graph $G''$ as follows. We delete from $G'$ the vertex $a$ and the vertices of $C_d$. We add the edge $vx$ and we add an edge from the neighbor $z_2'$ of $z_2$ outside $C_d$ to
$z_1'$. Then we apply induction where $x$ plays the role of $b_0$.

We may assume that there is no dangerous cycle in $G''$ containing the new edge $z_1'z_2'$. For then $G$ would contain a cycle $z_1'z_1z_2z_2'z_2''z_2'''w_1z_1'$ where $z_2',z_2'',z_2'''$ are joined to $C$. Then we form a new graph
 from $G'$
 by deleting $a,u,z_4,z_3,z_2,z_2',z_2''$ and adding the edge $vx$ and two edges from each of $z_2''',z_1$ to $C$. After the induction
 (where $x$ plays the role of $b_0$)
  we give $a$ the blue color of $x$, we make $z_2,u$ red and $z_4,z_3,z_2',z_2''$ blue. If $z_1'$ is blue, we change the color of $z_1$ to red. If $w_1$ is blue we change it to red. Now it is easy to color the blue vertices $1,2,3$: We first ignore the colors $1,2,3$ of $z_1,z_2'''$. Then we assign colors $1,2,3$ to $z_1$ (if it is blue),$z_4,z_3,z_2',z_2''',z_2''$ in that order or (if $z_1'$ is blue) $z_2''',z_2'',z_2',z_3,z_4$ in that order. So assume that there is no dangerous cycle.

By similar arguments, we may also assume that the new edge $z_1'z_2'$ is not part of a double edge. Also, $G''$ satisfies the connectivity condition $(c_1)$. So we can apply induction to $G''$.

After the induction $a$ will get the blue color of $x$, and we shall color $C_d$ as follows. $z_3,z_4$ are blue, and $u$ is red. If $z_1'$ is red, then $z_1$ is colored blue and $z_2$ is colored red. If $z_1',z_2'$ are blue, then then $z_1,z_2$ are colored red. If $z_1',z_2'$ are blue and red, respectively, then we make $z_1$ is blue and $z_2$ red unless we create a red facial $4$-path. In that case we make $z_2$ blue and $z_1$ red. Now it is easy to color the blue vertices of $C_d$ by the colors $1,2,3$ unless $z_1,z_2'$ are blue, and $z_1'$ is red.
(It is here useful to note that $w_2,w_3$ are blue, and only one more vertex of $C_d'$ can be blue.)
Below we dispose of that case.

In order to color $z_1,z_3,z_4$ we investigate the Kempe chains containing $z_2'$. Assume that $z_2'$ has blue color $1$.
We claim that the Kempe chain of colors $1,2$ (respectively $1,3$) containing $z_2'$ does not contain $b$ or $v$ or any blue neighbor of $v$. We show this by focussing at $C_d'$.
%Assume first that $z_1'$ is blue. Then $w_1,w_4$ are red. ($w_2,w_3$ are blue because they are joined to $C$.) If the third neighbor of $w_1$ is blue, it has the same color as $w_3$. If the third neighbor of $w_4$ is blue, it has the same color as $w_2$. Now it is clear that a Kempe chain cannot join the third neighbors of $w_1,w_4$.
 Recall that $z_1'$ is red. Then at most one of $w_1,w_4$ is blue
  because $w_1,w_2,w_3,w_4$ cannot all be blue as they are in a $5$-cycle.
  If $w_1$, say, is blue, then we may assume that the third neighbor $w_1'$ of $w_1$ is red since otherwise, we make $w_1$ red. There now may be a Kempe chain from $w_1$ to the third neighbor of $w_4$ which we denote $w_4'$. But then $w_1,w_4'$ have the same color, and we switch colors
of $w_2,w_3$ so that we prevent the Kempe chain under consideration to connect $w_1,w_4'$.

We shall now color $z_1,z_3,z_4$. We may assume that one of $w_1,w_4$ is blue since otherwise it is easy to color $z_1,z_3,z_4$ (possibly after switching colors in a Kempe chain so that $z_2',a$ have distinct colors). In the discussion below we assume it is $w_1$.
(Note that only one of $w_1,z_1',w_4$ can be blue because there cannot be four blue vertices in a $5$-cycle.)

If $z_2',w_1,a$ are colored $1,1,2$ or $1,2,2$, then we give $z_1,z_3,z_4$ colors $3,2,1$.

If $z_2',w_1,a$ are colored $1,2,3$, then we switch colors of the Kempe chain with colors $2,3$ containing $w_1$ so that $z_2',w_1,a$ are colored $1,3,3$, a previous case. (If $w_4$ is blue and $w_1$ is red we consider instead the Kempe chain with colors $1,2$ containing $z_2'$.)

If $z_2',w_1,a$ are colored $1,2,1$, then we switch colors of the Kempe chain with colors $1,3$ containing $z_2'$ so that $z_2',w_1,a$ are colored $3,2,1$, a previous case.

If $z_2',w_1,a$ are colored $1,1,1$, then we switch colors of a Kempe chain so that $z_2',w_1,a$ are not colored with the same color, a previous case.

\vspace{3mm}

{\bf Subcase $1.2.2.2.2.3.2.2$: $b$ is contained in $C_d'$.}

\vspace{3mm}

 By $(c_1)$, $b$ is a neighbor of $v$ on $C_d'$. Then $w_3=b,w_4=v$. We do induction as in Subcase $1.2.2.2.2.3.2.1$ except that we need not use Kempe chains. (Also, we may assume that we do not create a dangerous cycle
  containing the edge $z_1'z_2'$
  as $V(G)$ then consists of the vertices of $C$, the vertices $r_0',a$, the vertices of $C_d,C_d'$ and the vertices of the new dangerous cycle. It is then easy to color $G$ directly.) In Subcase $1.2.2.2.2.3.2.1$ we used Kempe chains only in the case where, after the induction, $a,z_1,z_2'$ are blue and $u,z_2,z_1'$ are red.
 $z_3,z_4,z_1$ do not yet have colors $1,2,3$. We ignore the colors $1,2,3$ of $a,w_2,w_3=b$. Possible $w_1$ or its third neighbor $w_1'$ is blue. But they are not both blue since, otherwise, we make $w_1$ red. Now we can assign colors $1,2,3$ to the vertices $z_1,z_3,z_4,a,w_3,w_2$ in that order.

% We recolor $z_1$ so that it has no color conflict, except possibly with $a$. If one of $w_4,z_1'$ is blue, we color $a$ red. Otherwise we keep $a$ blue but recolor $a$ so that it has a blue color distinct from the blue colors of $w_3=b$ and $z_1$. Now we can color $z_3,z_4$ unless $z_2',a$ have the same blue color, say $1$, and $z_1$ has color $2$, say. If $w_1$ is blue we make $z_1$ red and color $z_3,z_4$. So assume $w_1$ is red. Now we recolor $z_1,a$ so that $z_2',z_1,a$ have distinct colors, and now we can color $z_3,z_4$.

 \vspace{3mm}

This completes the proof of Case $1$.

\vspace{3mm}

Below we consider the Case $2$ where $b_0$ exists, and finally Case $3$ where $r_0$ exists and is $4$-forbidden. We already note here that Case $3$ is almost the same as Case $2$. In fact Case $3$ is much easier, because much of the reasoning in Case $2$ is about coloring the vertex $b_0$. This difficulty does not exist in Case $3$.

\vspace{3mm}

{\bf Case $2$: $b_0$ exists.}

\vspace{3mm}

 If all vertices in $int(C)$ are joined to $C-b_0$, then $G-V(C)$ is a path whose vertices are blue. It is part of a path starting at $b_0$. We traverse this path from $b_0$ and color its vertices $1,2,3,1,2,3, \ldots$.

 So assume that some vertex in $int(C)$ has no neighbor in $C$. Then there is a unique path $b_0b_1 \ldots b_kb_0'$ such that each of $b_0,b_1, \ldots ,b_k$ has a neighbor in $C-b_0$ but $b_0'$ has no neighbor in $C-b_0$. Possibly $k=0$. Let $H$ be the graph obtained from $G$ by deleting the vertices $b_1,b_2, \ldots ,b_k$ and adding the edge $b_0b_0'$. Let $x,y$ be the two other neighbors of $b_0'$ such that the path $b_0b_0'x$ turns sharp right at $b_0'$.

 We may assume that $G$ has no dangerous cycle. For such a cycle would have the form $C'=b_0b_1 \ldots b_kb_0'xb_1'b_2' \ldots b_m'ab_0$ where $b_0',a$ are the only two non-blue vertices.
 If $G-(V(C) \cup V(C'))$ is disconnected, then the two components are joined to $b_0',b_m'$, respectively. We apply induction to each component where $b_0'$ plays the role of $r_0$ and is $4$- forbidden, and $b_m'$ plays the role of $b_0$ (and a facial cycle of $C \cup C'$ plays the role of $C$). If $G-(V(C) \cup V(C'))$ is connected and is joined to only $b_0'$ in $C'$, then we again apply induction to $G-(V(C) \cup V(C'))$ where $b_0'$ plays the role of $r_0$ and is $4$-forbidden. Finally, if $m=0$ and $G-(V(C) \cup V(C'))$ is joined to each of $b_0',x$ by an edge, then the dangerous cycle $C'$ has the form $b_0b_1 \ldots b_kb_0'xab_0$. Since $int(C)$ has no
 $5$-path with five blue vertices, by Claim 1, it follows that
  $C'=b_0b_1b_0'xab_0$.
 (Note that a dangerous cycle cannot have $6$ or $7$ vertices. It may have $8$ vertices, but then it either has a blue $5$-path, which contradicts Claim 1, or else the two non-blue vertices have distance at least $3$ on the cycle.)
 Let $C''$ be the cycle obtained from $C \cup C'$ by deleting the edge $b_0a$. We apply, if possible, induction to $C''$ and its interior such that $x$ plays the role of $b_0$. If the induction is possible, it is easy to color $b_1,b_0$. So we may assume that induction is not possible, that is, $int(C''$) has a dangerous cycle $S$ containing the neighbor $y$ of $b_0'$ in $int(C'')$. If $y$ has a neighbor $y'$ on $S$ which is joined to $C$ (and therefore blue), we apply induction to the graph $G'$ obtained from $G-b_0'-b_1$ by adding an edge from $y$ to $b_0$ and an additional edge from $x$ to $C$. Then we apply induction so that $b_0$ plays the role of $r_0$ which is either right-forbidden or left-forbidden. We make the choice so that $y'$ is forced to be blue by that choice. However, as $y'$ is already blue because it is joined to $C$, we do not create a new dangerous cycle and therefore induction is possible. As $y$ becomes red, it has no color conflict with the blue vertex $x$. Now $b_0'$ is made red, and it is easy to color the blue vertices $b_1,b_0$. So we may assume that $y$ has no neighbor on $S$ which is joined to $C$. So $S=yz_1z_2z_3z_4y$ where $z_2,z_3$ are joined to $C$.
 Now we apply induction to each of the two components $H_1,H_2$ of $int(C)-b_1-b_0'-y-z_2-z_3$ where $H_1$ contains $z_1$ and $H_2$ contains $z_4,x$.
  (If $z_4,x$ are in distinct components, the argument is similar and easier.)
 Before induction we identify $z_i$ with a vertex on the outer cycle (and add an additional edge from $x$ to the outer cycle). We let $z_i$ play the role of $r_0$ and let it be $4$-forbidden for $i=1,4$. Then we make $y$ blue and $b_0'$ red. Finally we assign colors $1,2,3$ to $y,z_3,z_2,b_1,b_0$ in that order ignoring the neighbor of $z_1$ in $H_1$. If that neighbor has the same color as $y$ or $z_2$, then we permute colors in $H_1$. This shows that $G$ has no dangerous cycle.

\vspace{3mm}

{\bf Subcase $2.1$: $H-V(C)-b_0'$ is disconnected with components $M_x$ containing $x$ and $M_y$ containing $y$.}

\vspace{3mm}

  For $z \in \{x,y\}$, we let $H_z$ be the subgraph of $G$ induced by $V (C) \cup
V(M_z)$, together with the edge $zb_0$. We apply induction to each of the
graphs $H_x,H_y$. In each case $C$ is the outer cycle and $b_0$ is blue. If there
is one, let $x'$ be the only other blue vertex, other than either $b_1$ and $b_2$
or $x$ and $y$, that has distance at most $2$ from $b_0$ in $G$. Assume without loss of generality that $x'$ is in $M_x$. Apply the induction to $H_x$ and $H_y$. Suppose that $b_0'$
can be made red without violating Conclusion $(i)$. Notice that $x$ and $y$ and $x'$ are the only
possible neighbours of $b_k$ in the blue square-graph that are already
colored with one of the blue colors. We temporarily ignore the color of
$y$ and color in this order $b_k, \ldots ,b_1$ and recolor $b_0$ with the three blue
colors. There are at least two ways to do this; for one of these colorings,
$b_0$ has a color different from that of $x'$. If $y$ is blue, we complete the coloring by
permuting the blue colors in $M_y$ so that $y$ gets a color different from
those of $x$ and $b_k$.
If $x$ is red, then we color $b_0'$ with
the blue color of $b_0$ in $H_x$, and, again ignoring the color of $y$, color
$b_k, \ldots, b_1$ and recolor $b_0$. As in the preceding paragraph, one of the two
possibilities will give $b_0$ a color different from that of $x'$. If $y$ is blue,
appropriately permuting the blue colors of $M_y$ will give $y$ a blue color
different from those of $b_0$ and $b_k$, as required.
Thus, we may assume that:
\vspace{3mm}

{\bf Fact $y$}: $x$ is blue, $y$ is red, and making $b_0$ red introduces a red facial $4$-path $b_0 yz_1z_2$, with $z_1, z_2$ in $M_y$.

\vspace{3mm}

{\bf Subcase $2.1.1$: $M_y-y$ is disconnected with components $M_1,M_2$.}

\vspace{3mm}

Then we apply induction to $M_1,M_2$.
Before induction, we add an edge from the neighbor of $y$ in $M_1$ (respectively $M_2$) to $b_0$.
We get a red-blue coloring of $M_y$ with no facial red $4$-path, and, by permuting colors in $M_1$ we may assume that we get a proper coloring in the blue vertices in the square-graph of $M_y$ such that the color of $y$ is blue. However, this contradicts Fact $y$. So in Subcase $2.1$ we are left with

\vspace{3mm}

{\bf Subcase $2.1.2$: $M_y-y$ is connected.}

\vspace{3mm}

\vspace{3mm}

{\bf Subcase $2.1.2.1$: $y$ has a neighbor $y'$ in $M_y$ such that, if $y'$  is made blue, then we create no dangerous cycle in $M_y$.}

\vspace{3mm}

In this case we can apply induction to $M_y$ as follows. We let $y$ be joined to a vertex on the outer cycle which will play the role of $r_0$ and which will be either right-forbidden or left-forbidden depending on whether $y'$ is reached by a right turn or left turn on the path $r_0yy'$.
   Note that condition $(c_9)$ is satisfied because of the assumption of  Subcase $2.1.2$: $M_y-y$ is connected.
 After this induction $y$ is red. We now combine this coloring of $M_y$ with the coloring of $M_x$. Recall that $x$ is blue in this coloring
  by a remark immediately before Fact $y$.
 We let $b_0'$ be red which creates no red facial $4$-path in $M_y$. Then we color $b_k,b_{k-1}, \ldots ,b_0$ blue and we color them $1,2,3$ in that order.

So in Subcase $2.1.2$ and hence also in Subcase $2.1$ we are left with

\vspace{3mm}

{\bf Subcase $2.1.2.2$: If a neighbor $y'$ of $y$ in $M_y$  is made blue, we create a dangerous cycle $C'$. If the third neighbor of $y$, say $y''$, is made blue, we create a dangerous cycle $C''$ }

\vspace{3mm}

Assume first that $C'$ does not contain $y$.

We applied earlier induction to $M_y$ where $y$ is joined to the vertex playing the role of $b_0$. We concluded (Fact $y$) that $y$ would be red and that there would be a red facial $3$-path starting at $y$. Let us assume that in some such coloring the red facial $3$-path starting at $y$ turns sharp right. The case where it turns sharp left is treated in the same way. There may be another coloring where it turns sharp left, but we shall not use that. Let $y'$ be the neighbor of $y$ in this path.

We combine the two colorings obtained by applying induction to $M_x, M_y$ where each of $x,y$ is joined to a vertex playing the role of $b_0$. Then we make $b_0'$ red, and we obtain a facial $4$-path $b_0'yy'u$. Each of $y',u$ is in $C'$. There is at most one more vertex $v$ in $C'$ which is not joined to $C$ because $C'$ becomes dangerous if we make $y'$ blue. The vertex $v$ must exist because $M_y-y$ is connected
 (since otherwise, we would get a contradiction to Fact $y$).
Also, $v$ must be a neighbor of $y'$. Possibly $v$ is blue.

As $M_y-y$ is connected, it is not possible that both neighbors of $y''$ in $C''$ are joined to $C$. So, $y''$ has a neighbor $y'''$ on $C''$ which is not joined to $C$. If we make $y''$ blue, we make $C''$ dangerous. If we make $y'',y'''$ blue we make $C''$ forbidden. So one of $y'',y'''$ is red in the coloring of $M_y$. Hence, $y$ has at most two blue neighbors (in the square-graph) in $C''$. Because of Fact $y$, we conclude that $y$ must have a third blue neighbor in the square-graph. That blue vertex must be $v$ (because $y',u$ are red). And we shall get a problem if we make $v$ red (since otherwise, we could color $y$ blue, contradicting Fact $y$.) Summarizing, if we let $v,y$ interchange colors, then $y$ has no conflict with its blue neighbors (in the square-graph) in $C''$. Moreover, we must create a facial $4$-path containing $v$. This implies that $v$ has only one blue neighbor (in the square-graph) not in $C'$.
We now recolor as follows: We first give $y$ a blue color such that it has no color conflict with a blue vertex in $C''$. Then we keep $v$ as blue, but we color $v$ such that it has no color conflict with the blue vertices outside $C'$. We keep the blue color of the vertices of $C'$ joined to $C$ but we recolor them (using colors $1,2,3$) such that they have no color conflict with $v$ or the neighbor $u'$ of $u$ outside $C'$. This coloring contradicts Fact $y$.

In the argument above we assume that $C'$ does not contain $y$. So, we need to comment on the case where $C'$ contains $y$. In that case $C''=C'$, and $y,y',y''$ are the only vertices of $C'$ that are not joined to $C$. We apply induction to each of the two components of $H_y-V(C')$. We think of $y'$ (respectively $y''$) as a vertex on the outer cycle playing the role of $r_0$ and being $4$-forbidden. Then we can let $y$ be blue and obtain a contradiction to Fact $y$.

This completes the discussion of Subcase 2.1.

%. Now the only possible color conflict is that a vertex of $C'$ joined to $C$ has the same color as the neighbor $u'$ of $u$ outside $C'$. But then we permute the blue colors of the component of $M_y-V(C')$
%containing $u'$ so that $u'$ get a different blue color. This coloring contradicts Fact $y$ and completes the discussion of Subcase 2.1.

So, to complete the proof of Theorem \ref{t1} in Case $2$, there only remains Subcase $2.2$ below.

\vspace{3mm}

{\bf Subcase $2.2$: $H-V(C)-b_0'$ is connected.}

\vspace{3mm}

{\bf Subcase $2.2.1$: If $x$ is made blue, we create no dangerous cycle in $G-b_0$.}

\vspace{3mm}

We have already discussed the case where $r_0$ exists and is either right-forbidden or left-forbidden. We shall therefore apply Theorem \ref{t1} to $H$ where $b_0$ plays the role of $r_0$ and is right-forbidden. This is possible because of the assumptions in Subcases $2.2$ and $2.2.1$. After we we have applied Theorem \ref{t1}, $x$ is blue, and $b_0'$ is red. If also $y$ is red, or if $k> 0$, (or both), then we can color
$b_k,b_{k-1}, \ldots ,b_0$ blue, and we color them $1,2,3$ in that order. So, we may assume that $k=0$ (and hence $H=G$) and also the following:

\vspace{3mm}

\emph{Fact $x$ : If, in Subcase 2.2.1, we apply Theorem \ref{t1} to $G$ such that $b_0$ plays the role of $r_0$ (and it is either right-forbidden or left-forbidden), then $b_0'$ is red, and $x,y$ are blue. Moreover, there is a blue vertex $x'$ of distance $2$ to $b_0$, and $x',x,y$ have distinct blue colors.}

\vspace{3mm}

The existence of $x'$ in Fact $x$ combined with the assumption of Subcase $2.2$ implies that not both of $x,y$ are joined to $C$. If $y$ is joined to $C$, then we let $x,y$ interchange roles.
So, we may assume that $y$ is not joined to $C$.

Let $ab_0d$ be a $3$-path of $C$ traversed anticlockwise. We now form a new graph $H'$ from $G$ by deleting $b_0,b_0'$ and adding two new vertices $b,c$, and also adding the path $abcd$ and the edges $by,cx$.
We now try to use the induction hypothesis to $H'$ with $b$ playing the role of $r_0$ being left-forbidden or right-forbidden. (Note that $H'$ has more vertices than $G$ but fewer edges inside the outer cycle so it makes sense to try induction.) If this induction works, then $y$ will become red. We use this coloring to $G$ except that we let $b_0'$ be red and $b_0$ blue. Now give $b_0$ a color distinct from those of $x,x'$. This coloring would satisfy the conclusion of Theorem \ref{t1}, and therefore, we may assume that we cannot apply induction. This implies that we must have one of the following subcases $2.2.1.1$ or $2.2.1.2$ below.

\vspace{3mm}

{\bf Subcase $2.2.1.1$: $G-V(C)-y$ is disconnected.}

\vspace{3mm}

Note that then we do not create a dangerous cycle by making $y$ blue. This means that $x,y$ can interchange role unless $x$ is joined to $C$.
So, we may assume that $x'$ belongs to the component of $G-V(C)-y$ not containing $x$. (This is clear if $x$ is joined to $C$. And if $x$ is not joined to $C$, then we let $x,y$ interchange roles. If $G-V(C)-x$ is connected we proceed to Subcase $2.2.1.2$.)

Let us now focus on a coloring discussed in Fact $x$.
 (Now $x$ and $y$ have their original roles.)
 In this coloring $x$ and $y$ are blue, but $y$ cannot have two blue neighbors (in $G$). For if that were the case we could make $y$ red and we would have a contradiction to Fact $x$. Let the blue colors of $x,y,x'$ be $1,2,3$, respectively. We now interchange colors in the Kempe chain with colors $1,3$ in the square-graph containing $x$ that is, the connected component containing $x$ in the subgraph of $G^2$ consisting of the vertices of color $1,3$. Then $x$ changes color from $1$ to $3$. As $y$ has color $2$ and may have a neighbor (in $G$) of color $1$ or $3$, but not both, $x'$ keeps its color $3$. This contradicts the last statement in Fact $x$ and completes the discussion of Subcase $2.2.1.1$.

\vspace{3mm}

{\bf Subcase $2.2.1.2$: $G-V(C)-y$ is connected. If a neighbor of $y$ distinct from $b_0'$ is made blue, then we create a dangerous cycle in $H'$.}

\vspace{3mm}

Let the two other neighbors of $y$ be denoted $y',y''$. Let $C'$ be the dangerous cycle that arises if we make $y'$ blue. Let $C''$ be the dangerous cycle that arises if we make $y''$ blue.  Then all vertices of $C'$ are joined to $C$ except $y'$ and one or two more vertices. Again, we consider two subcases.

\vspace{3mm}

{\bf Subcase $2.2.1.2.1$: $C'=C''$.}

\vspace{3mm}

In this case $C'$ contains the non-blue vertices $y',y''$ and possibly a third non-blue vertex. Because of the connectivity conditions of $G$, that third non-blue vertex must exist and it must be $y$. So we may assume that $C'=yy'z_1z_2y''y$ traversed clockwise, and $int(C)-V(C')$ is disconnected with components $H_1$ joined to $y'$ and $H_2$ joined to $y''$
  (because $z_1,z_2$ are joined to $C$.)
Let us focus on a coloring discussed in Fact $x$. In this coloring $x$ and $y$ are blue. Also $z_1,z_2$ are blue. Then $y',y''$ are red. Let the colors of $x',x,y$ be $1,3,2$, respectively. In $G^2$ we consider the Kempe chain with colors $1,3$ containing $x'$. We may assume that this Kempe chain also contains $x$. If $x'$ is in $H_1$ this implies that the neighbors of $y',y''$ outside $C'$ have colors $1,3$, respectively.
 If $x'$ is in $H_2$, then we
  switch colors in the Kempe chain with colors $1,2$ containing $y$ so that $x',x,y$ get colors $1,3,1$, a contradiction unless that Kempe chain contains $x'$ which implies that the neighbor of $y''$ outside $C'$ has blue color $1$. In either case, the neighbor of $y''$ outside $C'$ is blue.
  %So assume that $x'$ is in $H_2$. Again, we consider the Kempe chain with colors $1,2$ containing $y$. We may assume that this Kempe chain also contains $x'$, that is, the neighbor of $y''$ outside $C'$ is blue and of color $1$.
  Now we interchange colors of $y,y'$. Then $y'$ is blue of color $2$. This may create a color conflict is with a vertex in $H_1$. However, then we recolor $H_1$ by applying the induction to $H_1$ with $y'$ playing the role of $b_0$.

\vspace{3mm}

{\bf Subcase $2.2.1.2.2$: $C' \neq C''$.}

\vspace{3mm}

In this case $C',C''$ are disjoint. We earlier introduced the facial path $ab_0d$ on $C$ traversed anticlockwise. One of $a,d$ is joined to $x'$ in $int(C)$. We may assume that $a$ is joined to $x'$. For if $d$ is joined to $x'$ in $int(C)$, then we interchange the roles of $x,y$. (This is possible because the existence of the disjoint cycles $C',C''$ imply that we do not create a new forbidden or dangerous cycle if we make $y$ blue.) We choose the notation such that $y'$ is obtained from a right turn at the edge $b_0'y$.
Now all vertices of $C'$ are joined to $C$ except $y'$ and one or two more vertices.

   Recall that $x$ is blue in $H'$. However,
it is not possible that $x$ is one of the blue vertices on $C'$. For if that were the case, then
   $x$ would be a neighbor of $y'$ by the connectivity condition $(c_1)$.
Then we let $x,y$ interchange roles and now we get a contradiction to the assumption of Subcase $2.2.1.2$
   because making $y'$ blue now does not create a new dangerous cycle.
   So, $x$ is not on $C'$. Clearly, $x'$ is not in $C'$ because of $C''$
    and the connectivity condition $(c_1)$.

As $G-V(C)-y$ is connected (since otherwise we are in Case $2.2.1.1$ which we have disposed of) it has a path from $x'$ to $x$. This path contains two vertices $z_1,z_2$ which are in $C'$ and distinct from $y'$ and the blue vertices on $C'$.

Let us again focus on a coloring discussed in Fact $x$.
%In this coloring $x$ and $y$ are blue, but $y$ cannot have two blue neighbors (in $G$).
Let the blue colors of $x,y,x'$ be $1,2,3$. respectively. One of the vertices $z_1,y',z_2$ may be blue, but two of them cannot be blue. In $G^2$ we consider the Kempe chain with colors $1,3$ containing $x'$. We may assume that this Kempe chain contains a path $P$ from $x'$ to $x$ since otherwise we can interchange colors $1,3$ in such a way that $x,x'$ have the same color which implies that we can color $b_0$ and complete
Subcase $2.2.1.2$ and hence Subcase $2.2.1$.
%We also consider the cycle $C''$ which becomes dangerous when we make the third neighbor $y''$ of $y$ blue. In this cycle $y''$ has two neighbors $u_1,u_2$ which are not joined to $C$, and all other vertices of $C''$ (except $y''$) are joined to $C$. (This follows by the same reasoning as we used for $C'$.) In the red-blue coloring we are focussing on, one of $y',y''$ is red. Let us assume that $y''$ is red. (If $y'$ is red we argue similarly.) Note that $x$ is distinct from $u_1,u_2$ by the assumption in Subcase $2.2.1$.

As only one of $z_1,z_2$ can be blue, the path $P$ cannot use the edge $z_1z_2$ in $G^2$. If $y'$ is blue,
 it has color $1$ or $3$, and
 then $P$ contains a path $z_1'y'z_2'$ where $z_1',z_2'$ are neighbors of $z_1,z_2$, respectively, not in $C'$. In that case we make $y'$ red. So assume $y'$ is red. Then  $P$ enters $C'$ in $z_1$ (or from the neighbor $z_1'$ of $z_1$ outside $C'$) and leaves $C'$ from $z_2$
(or goes to a neighbor $z_2'$ of $z_2$ outside $C'$). In between, $P$ uses only vertices of $C'$ which are joined to $C$. We may assume that only one of $z_i,z_i'$ is blue for $i=1,2$ since otherwise, we make $z_i$ red. Then we
  change
  colors of the blue colors of $C'$ joined to $C$. Now the Kempe chain with colors $3,1$ which contains $x'$ does not contain $x$ and hence $b_0$ can be colored. This completes the discussion of Subcase $2.2.1$.

 %Also, $G^2$ has a path $P'$ with colors $2,3$ from $y$ to $x$. This implies that one of $u_1,u_2$, say $u_2$ is blue and has color $3$. The cycle $C''$ has length $2$ modulo $3$, say length $5$. (Larger lengths are treated similarly or can be disposed of by a simple inductive argument.) Say $C''=y''u_1v_1v_2u_2y''$. Then $y_1,y_2$ are colored $1,2$ or $2,1$ because $u_2$ has color $3$. Because of $P$, a neighbor of $u_1$ has color $3$, and we also conclude that $y_1,y_2$ are colored $1,2$. We have now described the blue colors on $C''$ and it neighbors (except possibly a neighbor of $u_2$). But now we interchange  the way $C''$ and its neighbors are

\vspace{3mm}

{\bf Subcase $2.2.2$: If $x$ is made blue, we create a dangerous cycle $C_x$ in $G-b_0$. Similarly, if $y$ is made blue, we create a dangerous cycle $C_y$ in $G-b_0$.}

\vspace{3mm}

We may assume that $C_x,C_y$ are disjoint.
 For if $C_x,C_y$ share precisely one edge, then they have length $5$ and $H-V(C)$ consists only of $C_x,C_y,b_0'$ and two edges incident with $b_0'$, and it is easy to color $int(C)$. If $C_x,C_y$ share at least two edges, then $C_x=C_y$. If they do not contain $b_0'$, then $C_x=C_y= yxb_0''uvy$ where $u,v$ are joined to $C$. In this case we apply induction to $H-b_0'-x-y-v-u$ with $b_0''$ playing the role of $b_0$, and it is easy to extend the resulting coloring to $G$. Finally, if $C_x=C_y$, and $C_x$ contains $b_0'$,
we apply induction to the two components of $int(C)-V(C_x)$ such that each of $x,y$ plays the role of $r_0$ and is $4$-forbidden. It is then easy to color all the other vertices of $C_x$ and also the vertices $b_0',b_k, \ldots ,b_0$ blue and also color them by colors $1,2,3$.
 So we may assume that $C_x,C_y$ are disjoint.

Recall that $x'$ denotes the unique vertex in $int(C)$ which is a neighbor of $b_0$ in the square-graph. (We may assume that $x'$ exists since we otherwise replace an appropriate path on $C$ starting at $b_0$ by a single edge.) Since $H-V(C)-b_0'$ is connected, we may choose the notation such that, every path in $H-V(C)-b_0'$ from $x'$ to $x$ intersects $C_y$. (This includes the case where $x'$ is in $C_y$.) Let $y_1,y_2$ be the neighbors of $y$ in $C_y$ such that the path $b_0'yy_1$ turns left at $y$. Since $H-V(C)-b_0'$ is connected, $y_2$ is not joined to $C$. Possibly, $y_1=x'$.

\vspace{3mm}

{\bf Subcase $2.2.2.1$}: $y_1$ is not joined to $C$.

\vspace{3mm}

In this subcase $y,y_1,y_2$ are the only vertices of $C_y$ which are not joined to $C$. Hence $G-V(C)-V(C_y)$ has a component $Q_1$ containing $x'$ and a vertex $y_1'$ joined to $y_1$. Let $Q_2$ denote the other component of $G-V(C)-V(C_y)$ having a vertex $y_2'$ joined to $y_2$. We apply induction to each of $Q_1,Q_2$. When we apply induction to $Q_1,Q_2$ we add an edge from $y_i'$ to
a vertex of $C$ for $i=1,2$. With a slight abuse of notation we call that vertex $y_i$. When we apply induction, $y_i$ will play the role of $r_0$ and we call $y_i$ $4$-forbidden. We also add one more edge from $b_0'$ to $C$ in order to get a cubic graph. After the induction we make $y$ blue. Note that $b_0'$ is already blue and also has a color $1,2,3$. $y_1,y_2$ are now red. We then color $y$ and the vertices of $C_y$ joined to $C$ by the colors $1,2,3$. Finally, we color the vertices $b_k, \ldots ,b_0$ by the colors $1,2,3$. We may assume that precisely one of $b_0',x$ is blue since otherwise, we make $b_0'$ red. We now argue why this coloring of the blue vertices is possible.
We first color $y$. We ignore the blue colors in $Q_1$ at the moment. So, here we have to avoid the blue colors of $b_0',x,y_2'$. However, if both $b_0',x$ are blue, then we recolor $b_0'$ so that it becomes red. So we can color $y$. Then we color the other blue vertices of $C_y$ which is possible because we only have to watch $y_2',y$. Again, we ignore the blue colors in $Q_1$.
After this we permute the blue colors in $Q_1$ such that the color of $y_1'$ has no conflict with a blue color in $C_y$. Finally, we color the vertices $b_k, \ldots ,b_0$. There is only one problem, namely that we cannot color $b_0$ because of the color of $x'$ and two other vertices. We consider the case where these two other vertices are $y,x$ or $y,b_0'$. (In other words, we consider the case $k=0$. If $k\geq 1$, the argument is similar.) So, we may assume that $k=0$ and that $x,y,x'$ (or $b_0',y,x'$) have colors $1,2,3$, respectively. We may assume that $y_1'$ has color $3$ since otherwise, $y_1'$ has color $1$ (or is red), and we switch the two blue colors $2,3$ in $Q_1$. Then the neighbor $u_1$ of $y_1$ on $C_y$ joined to $C$ must have color $1$. The number of vertices of $C_y$ joined to $C$ is $2$ modulo $3$.
 By Claim 1,
we may assume that it is $2$ since otherwise, replace $4$ edges of $C_y$ with a single edge and use induction.
Again, we consider the Kempe chain with colors $1,3$ containing $x'$. As this must contain $x$ (or $b_0'$), and $y_2$ is red, we conclude that $y_2'$ is blue with color $1$. But then we can make $y$ red and color $b_0$ with the color $2$. (Note that the recoloring of $y$ from blue to red is because we want to give $b_0$ a color $1,2,3$. In Case $3$ below we shall encounter a similar situation with $r_0$ instead of $b_0$. Recoloring $y$ might give a red facial $4$-path. However, we need not recolor $y$ in Case $3$ below because we shall not give $r_0$ a color $1,2,3$.)

\vspace{3mm}

{\bf Subcase $2.2.2.2$}: $y_1$ is joined to $C$. That is, $y_1=x'$.

\vspace{3mm}

We may assume that $C_y=yy_2y_3y_4y_1y$
   (labelled anticlockwise)
 where $y_1$ is joined to $C$.

Thus we consider the case where $C_y$ has length $5$. (If $C_y$ has length $4$, that is, $C_y=yy_2y_3y_1y$ where $y_1,y_3$ are joined to $C$, then the proof is similar and easier. If $C_y$ is longer, then $C_y$ has three consecutive vertices joined to $C$. We replace these three vertices by a single edge and use induction. So only the case where $C_y$ has length $5$ needs consideration.)

By a similar argument, we may assume that $C_x=xx_2x_3x_4x_1x$
  (labelled clockwise)
   where two of $x_1,x_4,x_3$ are joined to $C$.

If $y_4$ is not joined to $C$, then $y_3$ is joined to $C$. Then we apply induction to the two components of $int(C)-y_4$. In one of the components we add the edge $y_1y_3$ before we use induction. In the other component we let $y_4$ play the role of $r_0$ and we call it $4$-forbidden. So we may assume that $y_4$ is joined to $C$. Now we delete from $H$ the vertices $y,y_1,y_4,y_3$ and use induction. Before the induction we also add an edge from the third neighbor $y_3'$ of $y_3$ to $C$. That neighbor on $C$ is also called $y_3$ with a slight abuse of notation. This $y_3$ will play the role of $r_0$ and is called $4$-forbidden. We add an edge from $y_2$ to $C$ and also an edge between $y_2,b_0'$. If the interior of $C$ is disconnected (after the deletion of $y,y_1,y_4,y_3$), then the neighbor $b_0$ of $b_0'$ on $C$ will play the role of $b_0$. (If the interior of $C$ is connected, which is the most difficult case, then $b_0$ is just a red vertex on $C$.) We do not create a dangerous cycle (because such a cycle would have to contain the edge between $y_2,b_0'$, and all vertices on this cycle outside $C_x$ would have to be joined to $C$.
Then $x_2,y_2$ have a common neighbor $y_2'$ joined to $C$. In this case it is easy to complete the proof by induction.)
We then apply induction. After the induction we make $y$ red. If $y_3'$ is red, it is easy to color the blue vertices $y_1,y_4,b_k, \ldots, b_0$ by $1,2,3$ in that order. (Before we do that we make $b_0'$ red in case $x$ is blue.) So assume that $y_3'$ is blue and has (possibly) the same color as $y_2$.

Now we try to make $y_2$ red. If we do not create a red facial $4$-path, then it is easy to color the blue vertices $b_k, \ldots, b_0,y_1,y_4$ by $1,2,3$ in that order. So assume that we create a red facial $4$-path if we make $y_2$ red. We keep $y_2$ blue, but give it a color $1,2,3$ distinct from its blue neighbors in $G^2$ except possibly $b_0'$. So now $y_2,b_0'$ may have the same blue color. If we can make $b_0'$ red we have finished. So assume that this is not possible, that is, either

(i) there is a red facial path $xx_1x_1'$ or else

(ii) the vertices $x,x_2$ are red.

In case (i) $x_3,x_4$ are joined to $C$. If we can make $x_2$ red, it is easy to color the blue vertices. So, $x_2'$ and one of its two other neighbors are red. But then ignore the colors $1,2,3$ of $x_3,x_4$ and we color the blue vertices $y_4,y_1,b_0', \ldots ,b_0,x_2,x_3,x_4$ in that order. So assume that we have (ii): $x,x_2$ are red.

Now we can give $b_0'$ a color $1,2,3$ distinct from that of $y_2$ and $x_1$. If possible, we make the color of $b_0'$ distinct from the color of $y_3'$. We color the blue vertices $y_4,y_1,b_k, \ldots, b_0$ by $1,2,3$ in that order. This works unless $y_3',y_2,x_1$ have distinct blue colors, say $1,2,3$, respectively. We uncolor $b_0'$ and try to switch colors in the Kempe chain with colors $2,3$ containing $x_1$. We may assume that then $y_2$ changes color. Then we try to switch colors in the Kempe chain with colors $1,3$ containing $x_1$. We may assume that then $y_3'$ changes color. But then $x_1,x_4,x_3,x_3'$ (where $x_3'$ is the third neighbor of $x_3$) are blue and of colors $3,2,1,3$ (or $3,1,2,3$), respectively. We try to make $x_3$ red. If this is possible, then we complete the proof using the afore-mentioned Kempe chain with colors $1,3$. So assume that the vertices in the facial $3$-path $x_2x_2'x_2''$ are red. We now make $x_2$ blue and $x_3,b_0'$ red. We ignore colors $3,2$ of $x_1,x_4$ and color now the blue vertices $x_2,x_4,x_1,y_4,y_1,b_0, \ldots, b_k$ in that order.

\vspace{3mm}

{\bf Case $3$: $r_0$ exists and is $4$-forbidden.}

\vspace{3mm}

This case is similar to and much easier than Case $2$. We use the same notation except that now we have a vertex $r_0$ instead of $b_0$.

 If $k>0$, then Case $3$ follows immediately from Case $2$ by first letting $r_0$ play the role of $b_0$ in Case $2$. After we obtain the desired coloring in Case $2$ we just change the color of $r_0$ from blue to red. So assume that $k=0$, that is, $H=G$.

The Subcase $3.1$ corresponding to Subcase $2.1$ is trivial: Just apply induction to $M_x$ and $M_y$ so that $b_0'$ becomes blue. In Subcase $3.2$ (where $H-V(C)-b_0'=G-V(C)-b_0'$ is connected) the Subcase $3.2.1$ is trivial
 because it reduces to Case 1.
So only Subcase $3.2.2$ needs attention. We repeat word for word the proof of Subcase $2.2.2$. In Subcase 2.2.2 we use an inductive argument in which $b_0'$ becomes blue. If $b_0'$ remains blue throughout the proof of Subcase $2.2.2$, then there is no problem in letting $r_0$ be red. In Subcase $2.2.2.1$ we make, at some stage, $b_0'$ red in order to be able to color $b_0$. This is not necessary in Subcase $3.2.2.1$. We also at some stage make $b_0'$ red if $x$ is blue. Also this is not a problem in Subcase $3.2.2.1$. The only case where an additional argument is needed is in Subcase $2.2.2.2$ just before statements (i),(ii). At this stage $y_3',y_2,b_0'$ are blue. $y_3',y_2$ have distinct blue colors, say $1,2$. The only color conflict is that $y_2,b_0'$ have the same color. We would like
   to change the color of $b_0'$ to either $3$ or red because then it is possible to color $y_4,y_1$ by colors $3,2$, respectively. So,
   we try to make $b_0'$ red. If this works, then we have finished in Subcase $2.2.2.2$. But in Subcase $3.2.2.2$ we now need to consider a possible red facial path $r_0b_0'xx_1$. As two vertices of $C_y$ are joined to $C$, we conclude that $x_3,x_4$ are joined to $C$. We may assume that $x_2$ is blue since otherwise we can give $y_4,
   y_1,
b_0'$ colors $3,1,3$, respectively. So, $x_2$ is blue of color $3$. If we can make $x_2$ red we have finished. So, we may assume that the third neighbor $x_2'$ of $x_2$ is red and that one of the two other neighbors of $x_2'$ is red. Now we can recolor $x_2,x_3,x_4$ with colors $1,2,3$ or $2,1,3$, respectively, so that the only possible color conflict is that $x_4$ has the same color as the third neighbor $x_1'$ of $x_1$. That color conflict can be eliminated by permuting colors $1,2,3$ in the component of $int(C)-x_1$ containing $x_1'$.

\section{Wegner's conjecture}

Theorem \ref{t1} is similar to Conjecture \ref{c1} except that the $3$-coloring of the blue graph is not obtained from Brooks' theorem. Also, Theorem \ref{t1} is very close to Wegner's conjecture when restricted to planar cubic $2$-connected graphs. Indeed, such a graph $G$ has a facial cycle of length at most $5$. We may assume that this is the outer cycle. We select a vertex on this cycle which we call $b_0$. We insert a vertex $d_0$ of degree $2$ on the outer cycle such that $b_0$ is
 adjacent to
the vertex $d_0$ of degree $2$ in order to satisfy condition $(c_4)$ in Theorem \ref{t1}. Then we apply Theorem \ref{t1}. The red square-graph is planar except for a pair of crossing edges in the outer face. The blue square-graph is $3$-colorable except that $b_0$ may have a blue neighbor in the blue square-graph when we ignore $d_0$. Thus, the square of $G$ can be colored in $7$ colors such that only two edges join vertices of the same color, a slight weakening of Wegner's conjecture.

To obtain the full version of Wegner's conjecture we need additional arguments. In this reasoning we shall use the classical result of Kotzig \cite{k} that every planar triangulation of minimum degree at least $4$ has a so-called \emph{light edge}, that is an edge such that the sum of degrees of its ends is at most $11$.

\begin{theorem}

\label{t2}

Let $G$ be a planar graph of maximum degree at most $3$. Then $G^2$ is $7$-colorable.

\end{theorem}

\emph{Proof of Theorem \ref{t2}}.
The proof is by induction on the number of vertices. The basis of the
induction is trivial so we proceed to the induction step. Assume (reductio ad absurdum) that Theorem \ref{t2} is false, and let $G$ be a counterexample with the smallest number of edges.
Clearly, $G$ has more than $7$ vertices.

\vspace{3mm}

\emph{Claim (1): $G$ is cubic and $2$-connected.}

\vspace{3mm}

\emph{Proof of Claim (1)}.
If $G$ has a vertex of degree $<3$, then it has degree at most $6$ in $G^2$.
 We delete this vertex. If the vertex has two neighbors and they are nonadjacent, we add an edge between its neighbors.
Then we use induction.
(The reason that we add an edge between the neighbors is that they are adjacent in $G^2$ and should therefore receive distinct colors when we use induction.) So, $G$ is cubic.

Clearly, $G$ is connected. If $G$ has a cut-edge $e$, then we delete $e$ and apply the
induction hypothesis to the connected components of the resulting
graph. By permuting the colors in one of the components, if necessary,
we obtain a $7$-coloring of $G^2$. So, $G$ is cubic and $2$-connected. This proves Claim (1).

\vspace{3mm}

\emph{Claim (2): $G$ is $3$-connected.}

\vspace{3mm}

\emph{Proof of Claim (2)}. Suppose (reductio ad absurdum) that $G$ contains two edges $x_1x_2,y_1y_2$ such that $G-x_1x_2-y_1y_2$ has two components $G_1,G_2$ such that $G_i$ contains $x_i,y_i$ for $i=1,2$. By choosing $G_1$ to be minimal we may assume that $x_1,y_1$ are not joined by an edge in $G$. If  $x_2,y_2$ are not joined by an edge in $G$, then we apply induction to (the square of) $G_i+x_iy_i$ for $i=1,2$. By permuting colors we may assume that $x_1,y_2$ have the same color, and $x_2,y_1$ have the same color. This results in a $7$-coloring of $G^2$, a contradiction. So assume that $x_2,y_2$ are joined by an edge in $G$, and let their third
neighbors
be $x_3,y_3$, respectively. Then we apply induction to (the square of) $G_1+x_1y_1$ and to $G_2-x_2-y_2+x_3y_3$. By permuting colors we may assume that $x_1,y_3$ have the same color, and $x_3,y_1$ have the same color.  By permuting the remaining colors we may assume that all colors of vertices adjacent (in $G$) to $y_3$ (except $x_3$) are also adjacent (in $G$) to $y_1$. Hence $y_2$ has $3$ available colors among the $7$ colors used for coloring $G^2$, and $x_2$ has at least one available color. Now we can color first $x_2$ and then $y_2$ and obtain a contradiction which proves Claim (2).

\vspace{3mm}

\emph{Claim (3): $G$ has no edge $xy$ which is contained in two distinct cycles $C_1,C_2$ such that $C_1$ has length $3$ and $C_2$ has length at most $5$.}

\vspace{3mm}

\emph{Proof of Claim (3)}.
Suppose (reductio ad absurdum) that $xy,C_1,C_2$ exist. As $G$ is $3$-connected, $C_1,C_2$ are facial cycles, and $C_2$ has length $4$ or $5$.
(Clearly, $C_1$ is facial, and clearly $C_2$ cannot have length $3$. If $C_2$ is nonfacial, then it has two vertices whose deletion makes the graph disconnected.)
If $C_2$ has length $4$ we contract $C_1,C_2$
into a vertex $v$ and use induction. If $v$ has color $1$ and the neighbors have colors $2,3,4$, then the vertices $x,y$ can receive two of the colors $2,3,4$, a third vertex of $C_2$ can be colored $1$, and now it is easy to color the two other vertices of $C_1 \cup C_2$ as well. So assume that $C_2=xx_1x_2x_3yx$ and $C_1=yy_1xy$. Delete the edge $xy$ and draw $G-xy$ such that the outer cycle is $C=xx_1x_2x_3yy_1x$. Let $x_1',x_2',x_3',y_1'$ be the neighbors of $x_1,x_2,x_3,y_1$, respectively, inside $C$. We may assume that $x_1',x_2',x_3',y_1'$ are distinct. For if two of them are identical, then we contract that vertex and $C$ into a single vertex (of degree $3$). We apply induction, and then it is easy to modify the coloring of the contracted graph to a $7$-coloring of $G^2$, a contradiction. Now we try to apply Theorem \ref{t1} to $G-xy$ where $y_1$ plays the role of $r_0$ and is right-forbidden. We also try to apply Theorem \ref{t1} to $G-xy$ where $y_1$ plays the role of $r_0$ and is left-forbidden.
If one of these attempts works, then we change the colors of $x,y$ to blue. As $y_1'$ is red, it is easy to give $x,y$ two colors $1,2,3$.
 So we may assume that it is not possible to apply Theorem \ref{t1}. Because $G$ is $3$-connected, $int(C)-y_1'$ is connected so $(c_6),(c_9)$ hold. The rest of $(c_1)-(c_9)$ all trivially hold, except $(c_8)$.
So we may assume that we create a dangerous or forbidden cycle $C'$ in $int(C)$
 when make a second neighbor $u_1$ of $y_1'$ blue.
 Possibly, $C'$ does not contain $u_1$.
Similarly we create a  dangerous or forbidden cycle $C''$ when make the third neighbor $u_2$ of $y_1'$ blue.
 Possibly, $C''$ does not contain $u_2$.
As a dangerous cycle has at least three blue vertices and there are only $4$ blue vertices when we apply induction, it follows that each of $C',C''$ has length $4$ or $5$. Hence they are facial cycles. As $C',C''$ have at least one of $x_1',x_2',x_3'$ in common, and there is only one facial cycle containing $x_i'$ and not intersecting $C$ it follows that $C'=C''$. As $G$ is $3$-connected, it follows that either $C'=C''=u_1u_2x_1'x_2'x_3'u_1$ in which case $G$ has $12$ vertices and $G^2$ has chromatic number $6$, or else $C'=C''=y_1'x_1'x_2'x_3'y_1'$ in which case $G$ has $10$ vertices and $G^2$ has chromatic number $6$,
 or else $C'=C''=x_1'x_2'x_3'z_1z_2x_1'$ where $z_1,z_2$ are distinct from $y_1',u_1,u_2$. In this case we contract all vertices of $C_1,C_2,C'$ into a single vertex and apply the induction hypothesis to the square of the resulting graph. The resulting $7$-coloring can easily be modified to a $7$-coloring
of $G^2$.

\vspace{3mm}

\emph{Claim (4): $G$ has no triangle.}

\vspace{3mm}

\emph{Proof of Claim (4)}.
Suppose (reductio ad absurdum) that $G$ has a triangle $x_1x_2x_3x_1$ which can be chosen to be the outer triangle. We now apply Theorem \ref{t1} where $x_1$ plays the role of $r_0$ and is $4$-forbidden. There is no dangerous cycle, as every dangerous cycle has at least three blue vertices. We may create a red facial path $x_2x_1x_1'x_2'$ or a red facial path $x_3x_1x_1'x_3'$
or both. In that case $x_2,x_2'$ and $x_3,x_3'$ are non-neighbors in the square-graph, by Claim (3). So, we identify $x_2,x_2'$ or $x_3,x_3'$ before we apply the $4$-Color Theorem to the red square-graph.
 %There are $6$ possibilities for this as $r_0$ may be left-forbidden or right-forbidden. When we have applied Theorem \ref{t1}, there may be a red facial path with $4$ vertices, but only one such path. When we identify the ends of this path we obtain a planar graph which is $4$-colorable unless it has a loop, that is, the ends of the red facial path with $4$ vertices are neighbors in the red square-graph. If this happens for each of the $6$ possible ways of using Theorem \ref{t1}, then the graph $G'=G-x_1-x_2-x_3$ has an outer cycle $C'$ with at most $6$ vertices. We apply Theorem \ref{t1} to $G'$ where $C'$ is red except that one vertex joined to the interior of $C'$ plays the role of $b_0$. (Note that the interior of $C'$ is connected because $G$ is $2$-connected.) After that we make all the vertices $x_1,x_2,x_3$ blue and we can also color them by the colors $1,2,3$ since at most two of them must avoid the color of $b_0$.
 This contradiction proves Claim (4).

\vspace{3mm}

\emph{Claim (5): $G$ has no non-facial cycle of length $<6$.}

\vspace{3mm}

\emph{Proof of Claim (5)}. Suppose (reductio ad absurdum) that $C$ is a non-facial cycle of length $<6$. By Claim (3), $C$ has no chord. Hence each edge not in $C$ but incident with a vertex of $C$ joins $C$ to a vertex inside or outside $C$. So precisely one or two edges join $C$ to its interior or exterior. This contradiction to Claim (2) proves Claim (5).

\vspace{3mm}

\emph{Claim (6): $G$ is cyclically $4$-edge-connected, that is, if $E$ is a set of three edges such that $G-E$ is disconnected, then $E$ consists of three edges incident with the same vertex.}

\vspace{3mm}

\emph{Proof of Claim (6)}.

Suppose (reductio ad absurdum) that $G$ has a set $E$ of three edges $x_1x_2, y_1y_2, z_1z_2$ such that $G-E$ has two components $G_1,G_2$ such that $G_i$ contains $x_i,y_i,z_i$ for $i=1,2$ and such that none of $G_1,G_2$ is a single vertex. By Claim $(4)$, each of $G_1,G_2$ has more
than
three vertices. We consider four new graphs $G_1',G_1'',G_2',G_2''$. $G_i'$ is obtained from $G_i$ by adding a vertex $g_i$ joined to $x_i,y_i,z_i$ for $i=1,2$.  $G_i''$ is obtained from $G_i$ by adding three vertices $x_i',y_i',z_i'$ forming a triangle and also adding the edges
$x_ix_i', y_iy_i', z_iz_i'$ for $i=1,2$. We first apply induction to $G_1''$ and $G_2''$. If $x_1,y_1,z_1$ get distinct colors in $G_1''$, and $x_2,y_2,z_2$ get distinct colors in $G_2''$, then it is easy to combine the two colorings to get a coloring of $G^2$, a contradiction. If $x_1,y_1,z_1$ get the same color in $G_1''$, then we apply induction to $G_2'$, and again it is easy to combine the two colorings. So we may assume that in $G_1''$ the vertices $x_1,y_1,z_1,x_1',y_1',z_1'$ have colors $1,1,2,3,4,5$, respectively. We may assume that in $G_2'$ the vertices $g_2,x_2,y_2,z_2$ have colors $1,3,4,5$, respectively. Now we try to combine the two colorings. The only possible conflict is that $z_1$ (which has color $2$ in $G_1$) can see a neighbor of $z_2$ in $G_2$ which also has color $2$. In $G_1$ we may switch colors $2,6$. We may also switch colors $2,7$. One of these two color switches will result in a proper coloring of $G^2$, a contradiction which proves Claim (6).

 %three $4$-cycles $C_1,C_2,C_3$ having a vertex $v$ in common. Then each $C_i$ has a unique vertex $x_i$ which is not in any of the other two $4$-cycles. Now we delete all vertices of $C_i-x_i$ for $i=1,2,3$, and we add edges so that $x_1x_2x_3x_1$ becomes a triangle $x_1x_2x_3x_4x_1$. Then we use induction. After the induction it is easy to extend the $7$-coloring so that we obtain a proper $7$-coloring of $G^2$.

\vspace{3mm}

  \emph{Claim (7): $G$ does not contain two distinct $4$-cycles having an edge in common.}

\vspace{3mm}

\emph{Proof of Claim (7)}. Suppose (reductio ad absurdum that $C_1:xx_1x_2yx$ and $C_2:xyx_3x_4x$ are $4$-cycles. We delete the edge $xy$ and think of $C_1\cup C_2-xy$ as a $4$-cycle. We apply induction to the resulting cubic graph. After the induction it is easy to color $x,y$. We may argue as follows: Let $x_i'$ be the third neighbor of $x_i$ for $i=1,2,3,4$. After the induction we may assume that $x_i$ has color $i$ for $i=1,2,3,4$. We can now give $x$ one of the colors $5,6,7$. Similarly for $y$. If the possible colors for $x,y$ are distinct, we have finished. So assume that $x_1',x_2',x_3',x_4'$ have colors $5,5,6,6$ or $5,6,5,6$, respectively. If we can change the color of $x_1$ to $3$, then we can give $x,y$ the colors $1,7$, respectively. So assume that the color $3$ is present at a neighbor of $x_1'$. Similarly, a neighbor of $x_2'$ (respectively $x_3'$, respectively $x_4'$) has color $4$ (respectively $1$, respectively $2$).
If we can change the colors of $x_1,x_3$ to $7$, then we can give $x,y$ the colors $1,3$, respectively. So, we may assume that the color $7$ is present at a neighbor of $x_1'$ and also at a neighbor of either $x_2'$ or $x_4'$ or both. Assume that $7$ is present at a neighbor of $x_2'$. (If $7$ is present at a neighbor of $x_4'$, the proof is similar.) We switch colors of $x_1,x_2$. We may assume that the color $2$ is present at a neighbor of $x_3'$ since otherwise, we could change the color of $x_3$ to $2$ and complete the proof. Similarly the color $1$ is present at a neighbor of $x_4'$. Now the vertices $x_1,x_2,x_3,x_4,x,y$ are colored $4,3,7,3,2,1$, respectively.

\vspace{3mm}

The dual version of Kotzig's result on light edges in triangulations implies that $G$ has two facial cycles $C_1,C_2$ of length $k_1,k_2$ respectively, such that $C_1,C_2$ have an edge $xy$ in common and such that $k_1 \leq k_2, k_1+k_2 \leq 11$. Hence $k_1 \leq 5$.
 By Claim (4), $k_1 \geq 4$. By Claim (7), $k_2 \geq 5$.

We choose $C_1,C_2$ such that $k_1+k_2$ is minimum. We delete the edge $xy$ and draw $G$ such that the outer cycle $C$ is $C_1 \cup C_2-xy$. This cycle can be described as $C:xx_1x_2 \ldots x_{k_2-2}yy_1y_2 \ldots y_{k_1-2}x$. Let the third neighbors of $x_1, x_2, \ldots ,y_{k_1-2}$ be denoted $x_1', x_2', \ldots ,y_{k_1-2}'$, respectively. As $G$ is $3$-connected, and there are at most $7$ edges from $C$ to its interior (since $k_1+k_2 \leq 11$), it follows from Claim (6) that
 $G-V(C)$ is connected.

 %has at most two connected components, and each of these is attached to both $C_1$ and $C_2$ and also to
%one of $x_1,x_{k_2-2}$ and also to one of $y_1,y_{k_1-2}$. So either $G-V(C)$ has one component $Q_0$ or else it has two components $Q_1,Q_2$.
%If it has two components, then we may choose the notation such that $Q_1$ contains $x_1',y_{k_1-2}'$, and $Q_2$ contains $y_1',x_{k_2-2}'$.

We now apply Theorem \ref{t1} to $G-V(C)$  where one of $x_1,x_{k_2-2},y_1,y_{k_1-2}$ plays the role of $r_0$ and all other vertices of $C$ are red. We call $r_0$ either left-forbidden or right-forbidden in order to prevent that there is a red facial $4$-path containing an edge of $C$ and starting at the neighbor of $r_0$ on $C$ distinct from $x,y$. (We shall later make $x,y$ blue so that a red facial path cannot start at $x$ or $y$.) We divide the argument into two cases.

Consider first the case where $k_2 \leq 6$. We apply Theorem \ref{t1} to $G-xy$ where we let $x_1$ (or $x_{k_2-2}$ or $y_1$ or $y_{k_1-2}$) play the role of $r_0$.
 Before we show that we can apply Theorem \ref{t1}, we explain how this will complete the proof. After the application of Theorem \ref{t1}
 we make $x,y$ blue, and we can extend the $3$-coloring of the blue square-graph to first $y$ and then $x$ because $x$ is adjacent (in $G^2)$ to at most one blue vertex inside $C$. To justify the last statement, the statement $(iii)$ in Theorem \ref{t1} implies that the neighbor of $x_1$ inside $C$ is red, and therefore it is possible to give $x$ a blue color and also a color $1,2,3$. We then apply the $4$-Color Theorem to the red square-graph. The only problem is that there may be a facial $4$-path when $k_2=6$, namely $x_1x_2x_3x_4$. Note that the vertices $x_1,x_4$ are not neighbors in the square-graph because of Claim (5). So, before we apply the $4$-Color Theorem we identify $x_1,x_4$. After this identification the red square-graph is planar.

We now explain why we can apply Theorem \ref{t1} to $G-xy$. Claims (4),(5),(6) and the minimality of $k_1+k_2$ imply that the vertices $x_1', x_2', \ldots ,y_{k_1-2}'$ are distinct.
We claim that condition $(c_9)$ is satisfied. That is, $int(C)-x_1'$ is connected. For suppose that $int(C)-x_1'$ has two components $H_1,H_2$. As $G$ is cyclically $4$-edge-connected, $k_1+k_2=11$, that is $k_2=6,k_1=5$, and the notation can be chosen such that $H_1$ contains $y_3',y_2',y_1'$, and $H_2$ contains $x_2',x_3',x_4'$. But then the edges $yy_1,xy_3$ and the edge from $x_1'$ to $H_1$ separate $G$, a contradiction to Claim (6). Similarly, $int(C)-x'$ is connected whenever $x$ is one of $x_1,x_{k_2-2},y_1,y_{k_1-2}$.

So, the only problem in the case $k_2 \leq 6$ is that there may be a forbidden or dangerous cycle $C_1'$ when we try to apply Theorem \ref{t1} with $x_1$ (or one of $x_{k_2-2},y_1,y_{k_1-2}$) playing the role of $r_0$ and being right-forbidden or left-forbidden. Then $C_1'$ is disjoint from $C$, and $C_1'$ contains at least three blue vertices. At least two of these are joined to $C$. We also try to apply Theorem \ref{t1} with  one of $x_{k_2-2},y_1,y_{k_1-2}$ playing the role of $r_0$ (and being right-forbidden or left-forbidden). Assume $C_2',C_3',C_4'$ are the resulting forbidden or dangerous cycles. Then each $C_i'$ is disjoint from $C$, and $C_i'$ contains at least three blue vertices. At least two of these are joined to $C$.

  We claim that precisely two vertices of $C_i'$ are joined to $C$. For, if three vertices of $C_1'$, say, are joined to $C$, then $C_1'$ has length $5$ and shares an edge with a $4$-cycle. Hence $k_1=4, k_2 \leq 5$. As we have noted earlier that $k_2 \geq 5$, the notation can be chosen such that $C_1'=x_2'x_3'y_1'z_1z_2x_2'$. But then there is no dangerous cycle when we think of $x_3$ as $r_0$ and let it be left-forbidden. This contradiction shows that $C_i'$ is joined to precisely two vertices of $C$ for each $i=1,2,3,4$.

Consider the graph $H$ induced by $C$ and $C_1'$ and $x_1'$. As $C_1'$ has at most two non-blue vertices, there are at most $7$ edges from $H$ to vertices not in $H$. As $G$ is cyclically $4$-edge-connected, at most one face of $H$ is not a face of $G$. Now it is easy to see that some facial cycle of $H$ is a $4$-cycle and hence $k_1 \leq 4$. If $C_1'$ does not contain $x_1'$, then $C_1'$ contains a path $x_1''x_2'x_3'$ where $x_1''$ is a neighbor of $x_1'$. Then the edge $x_2x_2'$ is contained in a facial $4$-cycle and a facial $5$-cycle implying that $k_2 \leq 5$. But now we can apply Theorem \ref{t1} to $G-xy$ with $x_3$ playing the role of $r_0$ and being right-forbidden or left-forbidden (because now there cannot be a dangerous cycle). So we may assume that $C_1'$ contains $x_1'$.
  But then at least three vertices of $C_1'$ are joined to $C$, a contradiction to an earlier claim.

%Similarly each of $C_2',C_3',C_4'$ are facial cycles containing three neighbors of $C$. Then all of $C_1',C_2',C_3',C_4'$ must be the same facial cycle. For if two of
 %{\color{red}them}
%are distinct, then $k_2=6$ and one must contain the vertices $x_1',x_2',y_2'$ and the other must contain $x_3',x_4',y_1'$. But then the edge $y_2y_2'$ is contained in a facial $5$-cycle %and hence $k_2 \leq 5$, a contradiction. So $C_1'=C_2'=C_3'=C_4'$.
%{\color{red}But then $C_1$ is joined to at least $3$ vertices of $C$, a contradiction to an earlier claim.}

%and it contains the path $x_1'x_2'x_3'x_4'$ (if $k_2=6$). As it also contains either $y_1'$ or a neighbor of $y_1'$, and it contains either $y_2'$ or a neighbor of $y_2'$, it follows that $C_1'$ is the cycle $x_1' \ldots x_{k_2-2}'y_1' \ldots y_{k_1-2}'x_1'$ showing that $k_1=k_2=4$. Now $G$ has only $10$ vertices and is easy to color $G^2$.

This completes the case $k_2<7$.

\vspace{3mm}

Consider finally the case where $k_2=7$. Then $k_1=4$.

%Consider first the subcase where $G-V(C)$ has two components $Q_1,Q_2$. As $G$ is $3$-connected, the notation can be chosen such that $Q_1$ contains $y_2',x_1',x_2',x_3'$, and $Q_2$ contains $x_4',x_5',y_1'$. We apply Theorem \ref{t1} to $Q_1,Q_2$ where $x_1,x_5$ play the role of $b_0$. It is clear that we can apply Theorem \ref{t1} to $Q_2$. If we cannot apply Theorem \ref{t1} to $Q_1$, then $Q_1$ contains a cycle containing $y_2',x_2',x_3'$ and two more vertices. But then we apply Theorem \ref{t1} to $Q_2$ where $x_1$ is either right-forbidden or left-forbidden.
%We make $x$ blue and $y$ red. Note that $x_5$ is blue. If $x_1,x_1'$ are both blue we change $x_1$ from blue to red. Then we can give $x$ a color $1,2,3$ because we can permute the blue colors in $Q_2$. If $x_1$ is red, we identify $x_1,x_4$ before we apply the $4$-Color Theorem to the red square-graph. So assume that $G-V(C)$ is connected.

We try to apply Theorem \ref{t1} where one of $x_1,x_5$ plays the role $b_0$. If this is not possible, then $G-V(C)$ contains a cycle which contains three of the vertices $x_2',x_3',x_4',y_1',y_2'$ and one or two more vertices. This is easily seen to contradict the assumption $k_2=7$ and that $G$ is $3$-connected and $G-V(C)$ is connected. (Consider for example the case where the dangerous cycle contains $y_2',x_2',x_3'$. Then the dangerous cycle must be of the form $y_2'r_1r_2x_2'x_3'y_2'$ because $G$ is $3$-connected. But then $G-V(C)$ is disconnected, a contradiction.) So we assume that we can apply Theorem \ref{t1} where $x_1$ plays the role $b_0$.

If the neighbor $x_1'$ of $x_1$ inside $C$ is red, then we make $x$ blue and we can give it a color $1,2,3$. Then we focus on the red square-graph. We first delete $y$ and then identify
$x_2,x_5,y_2$. Then we apply the $4$-Color Theorem. We can extend the $4$-coloring to include $y$ because $y$ can see only the colors of $x_4,x_5,y_1$.

If the neighbor $x_1'$ of $x_1$ inside $C$ is blue, then we make $x,y$ red. Then we focus on the red square-graph. We first delete $y,x$ and then identify
$x_2,x_5,y_2$ as before. Then we apply the $4$-Color Theorem. We can extend the $4$-coloring to include $y,x$ because $y$ can see only the colors of $x_4,x_5,y_1$, and $x$ can see only the colors of $y,y_1,y_2$.

This contradiction completes the proof of the theorem.

\vspace{3mm}

  {\bf Acknowledgements.} Thanks are due to the referees for numerous comments that greatly improved the presentation.


\begin{thebibliography}{99}



\bibitem{hjt} S.G.~Hartke, S.~Jahanbekam, and B.~Thomas,
The chromatic number of the square of subcubic planar graphs. Manuscript, April 2016.



\bibitem{bz}
Y.~Bu and X.~Zhu, An optimal square coloring of planar graphs.
{\em J.Combinatorial Optimization}
{\bf 24}
(2012)
580--592.

\bibitem{cj}
D.W.~Cranston and
B.~Jaeger, List-coloring the Squares of Planar Graphs
without 4-Cycles and 5-Cycles.
{\em arXiv:1505.03197v1 [math.CO]} 13 May 2015.


\bibitem{dkns}
Z.~Dvo\v r\'ak, D.~Kr\'al, P.~Nejedl\'y, and R.~\v Skrekovski, Coloring squares of planar graphs with girth six.
{\em European Journal of Combinatorics}
{\bf 29}
(2008)
838--849.


\bibitem{g}
M.~Gionfriddo, A short survey of some generalized colorings of graphs.
{\em Ars Combinatoria}
{\bf 30}
(1986)
275--284.


\bibitem{jt}
T.~Jensen and B.~Toft, Graph Coloring Problems.
{\em John Wiley, New York}
{1995}.


\bibitem{kp}
S.-J.~Kim and B.~Park, Coloring the squares of graphs whose maximum
average degrees are less than 4.
{\em arXiv:1506.04401v1 [math.CO]} 14 Jun 2015.

\bibitem{k}
A.~Kotzig, Contribution to the theory of Eulerian polyhedra.
{\em Mat.-Fyz. CP asopis. Slovensk. Akad. Vied} {\bf 5} (1955) 101-113.


\bibitem{mt}
B.~Mohar and C.~Thomassen,
Graphs on Surfaces.
{\em Johns Hopkins University Press, Baltimore}
(2001).


\bibitem{ms}
M.~Molloy and M.R.~Salavatipour, A bound on the chromatic number of the square of a planar graph.
{\em J. Combinatorial Theory Ser. B}
{\bf 94}
(2005)
189--213.


\bibitem{w}
G.~Wegner, Graphs with given diameter and a coloring problem
{\em preprint,University of Dortmund}
(1977).



\end{thebibliography}
\end{document}